\documentclass[12pt,reqno]{amsart} 
\usepackage{amsthm, amssymb,
  amsmath,verbatim,float,fullpage,mathpazo,psfrag,graphicx} 
\usepackage[dvipsnames]{xcolor}
\numberwithin{equation}{section}
\newtheorem{Theorem}{Theorem}[section]
\newtheorem*{Theorem*}{The Main Theorem} 
\newtheorem{Lemma}[Theorem]{Lemma}
\newtheorem{Proposition}[Theorem]{Proposition}

\newtheorem{Remark}[Theorem]{Remark}
 
\renewcommand{\le}{\leqslant} 
\renewcommand{\ge}{\geqslant} 
\newcommand{\lra}{\longrightarrow} 
\newcommand{\cbrac}{)\hspace{-1.8mm}(} 
\newcommand{\ux}{\mathbf x}
\newcommand{\exprS}{\mathcal S}
\newcommand{\exprT}{\mathcal T}
\renewcommand{\P}{\mathbf P}
\newcommand{\QQ}{\mathbf Q} 
\newcommand{\ZZ}{\mathbf Z}
\newcommand{\complex}{\mathbf C} 
\newcommand{\FF}{\mathcal F} 
\newcommand{\unord}{\mathcal Z} 
\newcommand{\cover}{\mathcal Y} 
\newcommand{\Gal}{\text{Gal}}
\newcommand{\SG}{\mathfrak S}
\newcommand{\conic}{\mathcal K}
 
\newcommand{\basefield}{\mathcal E} 
\newcommand{\pascal}[6]{\left\{ \begin{array}{ccc} #1 & #2 & #3\\ #4 &
     #5 & #6 \end{array} \right\}}
\newcommand{\pascarray}[6]{\left[ \begin{array}{ccc} #1 & #2 & #3\\ #4 &
     #5 & #6 \end{array} \right]}
\begin{document} 

\title[Reconstruction problem for Pascals]{On the reconstruction problem for Pascal lines} 
\maketitle 
\centerline{Abdelmalek Abdesselam and Jaydeep Chipalkatti} 

\bigskip 

\bigskip 

\parbox{17cm}{ \small
{\sc Abstract:} Given a sextuple of distinct points $A, B, C, D, E, F$ on a
conic, arranged into an array $\pascarray{A}{B}{C}{F}{E}{D}$, 
Pascal's theorem says that the points $AE \cap BF, BD \cap CE, AD \cap CF$ are collinear. The line
containing them is called the Pascal of the array, and one gets 
altogether sixty such lines by permuting the points. In this paper we
prove that the initial sextuple can be explicitly reconstructed from four
specifically chosen Pascals. The reconstruction formulae are encoded
by some transvectant identities which are proved using the graphical calculus for binary forms.} 

\bigskip


\medskip 


\bigskip 

\tableofcontents

\section{Introduction} 

This paper solves a reconstruction problem which arises in the context
of Pascal's hexagram in classical projective geometry. The main
result will be explained below once the required notation is
available. 
\subsection{} 
Let $\P^2$ denote the complex projective plane, and fix a nonsingular
conic $\conic$ in $\P^2$. Suppose that we are given six distinct
points $A, B, C, D, E, F$ on $\conic$, arranged as an array 
$\left[ \begin{array}{ccc} A & B & C \\ F & E & D \end{array}
\right]$. Then Pascal's theorem\footnote{One can find a proof in virtually any book on
  elementary projective geometry, e.g., Pedoe~\cite[Ch.~IX]{PedoeC} or
  Seidenberg~\cite[Ch.~6]{Seidenberg}. It is doubtful whether Pascal
  himself had a proof.} says that the three cross-hair intersection points 
\[ AE \cap BF, \quad BD \cap CE, \quad AD \cap CF\] 
(corresponding to the three minors of the array) are collinear. 

\begin{figure}
\includegraphics[width=8cm]{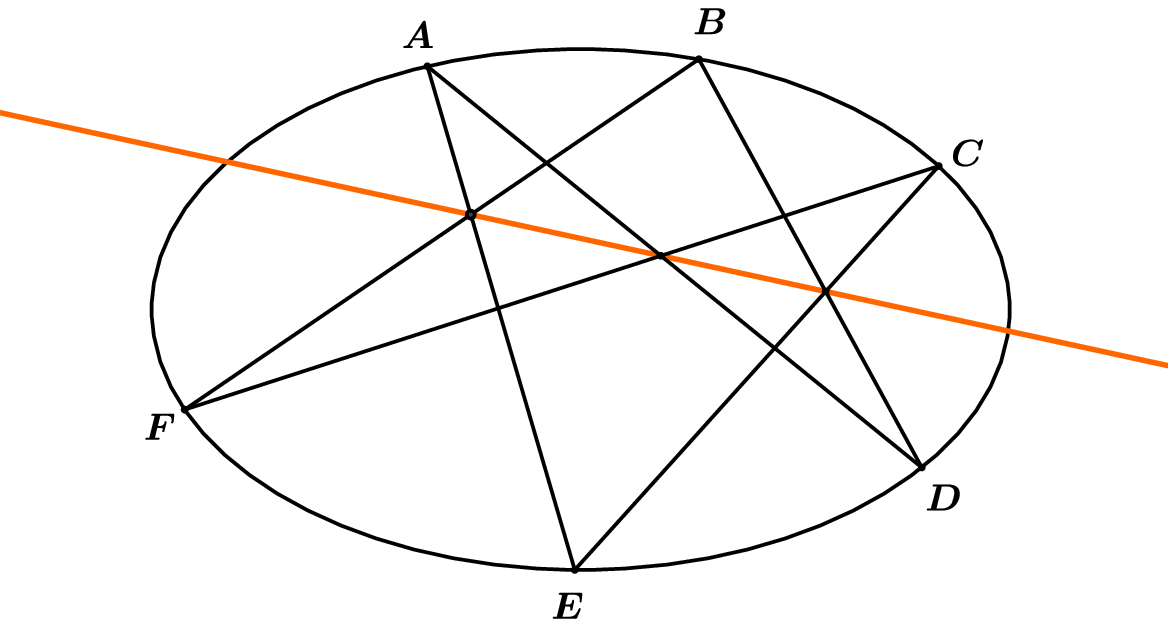}
\caption{Pascal's theorem} 
\end{figure}  

The line containing them is called the Pascal line, or just the Pascal, of
the array;  we will denote it by $\pascal{A}{B}{C}{F}{E}{D}$. It is
easy to see that the Pascal remains unchanged if we permute the rows
or the columns of the array; thus 
\begin{equation} 
\pascal{A}{B}{C}{F}{E}{D}, \quad \pascal{F}{E}{D}{A}{B}{C}, \quad 
\pascal{E}{D}{F}{B}{C}{A} 
\label{3pascals} \end{equation} 
all denote the same line. 

Any essentially different arrangement of the same points, say 
$\pascal{E}{A}{C}{B}{F}{D}$, corresponds \emph{a priori} to a different
line. Hence we have a total of $\frac{6!}{2!3!} = 60$ notionally
distinct Pascals. It is a theorem due to Pedoe~\cite{Pedoe}, that
these $60$ lines are distinct if the initial six points are
chosen generally.\footnote{If one tries to draw a diagram of 
  the sextuple together with all sixty of its Pascals, a dense and 
incomprehensible profusion of ink is the usual outcome. The curious reader is referred to 
{\tt http://mathworld.wolfram.com/PascalLines.html}} The configuration
of six points with all of its associated lines is sometimes called Pascal's
hexagram. 

The best classical references for the geometry of Pascal lines are
by Salmon~\cite[Notes]{Salmon} and Baker~\cite[Note II,
pp.~219--236]{Baker}. An engaging recent account is given in
the article by Conway and Ryba~\cite{ConwayRyba}. The reader is
referred to~\cite{KadisonKromann} and~\cite{Seidenberg} for standard
facts about projective planes. 

\subsection{} It is natural to wonder to what extent the construction sequence 
\[ \text{six points on $\conic$} \leadsto \text{sixty lines
  in the plane} \] 
can be reversed; that is to say, whether one can \emph{reconstruct} the
initial sextuple if the positions of some of the Pascals are
known.\footnote{The conic itself is fixed throughout, and as such assumed to
  be known.} In this paper we establish the following result: 

\begin{Theorem*}[Preliminary Form] \rm 
The sextuple $A, \dots, F$ can be reconstructed from the
following four Pascals: 
\begin{equation} 
\ell_1 = \pascal{A}{D}{B}{E}{C}{F}, \quad  
\ell_2 = \pascal{A}{C}{F}{E}{D}{B}, \quad  
\ell_3 = \pascal{A}{D}{F}{E}{C}{B}, \quad 
\ell_\star = \pascal{A}{B}{C}{F}{D}{E}. \label{list.fourpascals} \end{equation} 
\end{Theorem*} 
The arrays follow a pattern and the last one is on a
different footing from the first three; this will be explained in section~\ref{section.proof.overview}. 

\subsection{} \label{section.maintheorem.refined} 
In order to state the theorem more precisely, let $[z_0,
z_1, z_2]$ be the homogeneous coordinates on $\P^2$, and let the conic
$\conic$ be defined by the equation $z_1^2 = z_0 \, z_2$. Lines in
$\P^2$ are also given by homogeneous coordinates; for instance, 
the line $2 \, z_0 + 3 \, z_1 + 5 \, z_2=0$ has line coordinates
$\langle 2, 3, 5 \rangle$. 

Choose independent variables $a, \dots, f$, and fix the points 
\begin{equation} 
A = [1,a,a^2], \quad B =[1,b,b^2], \quad \dots, F = [1,f,f^2] 
\label{point.list} \end{equation} 
on $\conic$. 

Let $\langle 1, s_i, t_i \rangle$ denote the line coordinates
of $\ell_i$ for $i=1,2,3$, and $\langle 1, s_*, t_* \rangle$ those of $\ell_*$. 
Each of these Pascals is obtained by starting from the
points in (\ref{point.list}) and taking joins and intersections, hence
it is intuitively clear that $s_i$ and $t_i$ are rational functions in $a,
\dots, f$. The actual expressions are rather cumbersome; for
instance, 
\begin{equation} 
s_1 = \frac{abf - abe  - acd + \text{$9$ similar terms}}{abce-abcf 
  + \text{$4$ similar terms}},  
\quad 
t_1 = \frac{ac - af -bc + \text{$3$ similar terms}}{abce-abcf 
  + \text{$4$ similar terms}},  
\label{expressions.st} \end{equation} 
and likewise for the other $s_i, t_i$. The reconstruction problem is
to go backwards from the collection of Pascals $\{\ell_1, \ell_2, \ell_3,
\ell_*\}$ to the collection of points $\{A, \dots, F\}$. Our result
says that this can be done in algebraically the simplest possible
way. 

\medskip 

\begin{Theorem*}[Refined Form] \rm Each of the variables $a, \dots, f$ can be
expressed as a rational function of $s_i$ and $t_i$ for $i=1,2,3,*$. 
\end{Theorem*}

A naive attempt to prove the theorem would start from the
formulae for $s_1, \dots, t_*$, and try to `solve' for the variables $a,
\dots, f$. However, the expressions in~(\ref{expressions.st}) are too
complicated for this to succeed. We will instead use binary quadratic forms to represent points and lines in
$\P^2$, and express their joins and intersections in the language of
transvectants (see section~\ref{section.binary.forms}). One can then make these rational functions
completely explicit by exploiting the geometry of the Pascals in
conjunction with the graphical calculus for binary forms. It is an immediate corollary of the main theorem
that the Galois group of Pascal lines is isomorphic to the symmetric
group $\SG_6$. 

Our main theorem is thematically similar to, and partly inspired by, Wernick's problems in Euclidean triangle 
geometry - more on this in section~\ref{section.Wernick} below. 

\subsection{An overview of the proof} \label{section.proof.overview} 
\begin{figure}
\includegraphics[width=10cm]{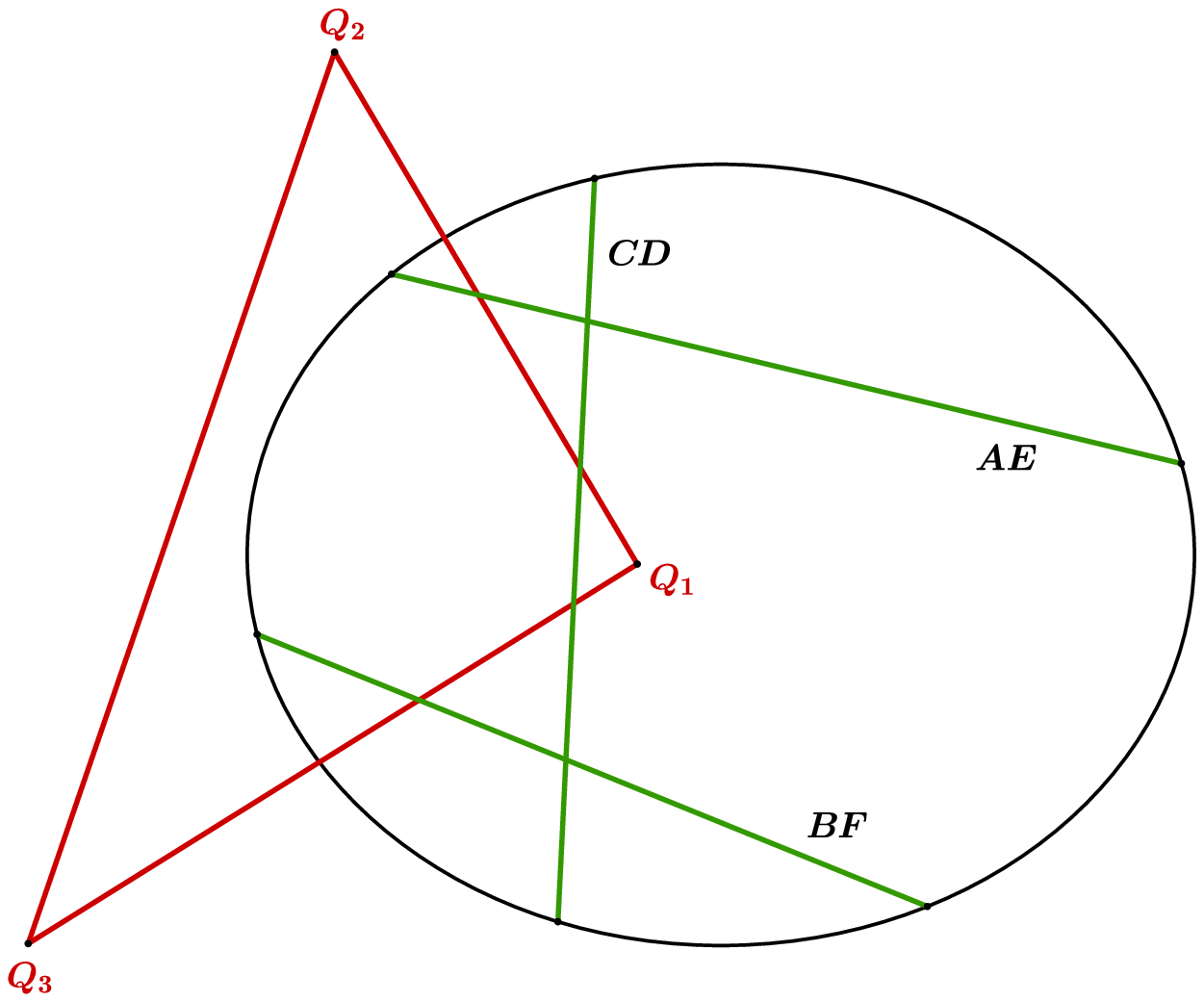}
\includegraphics[width=10cm]{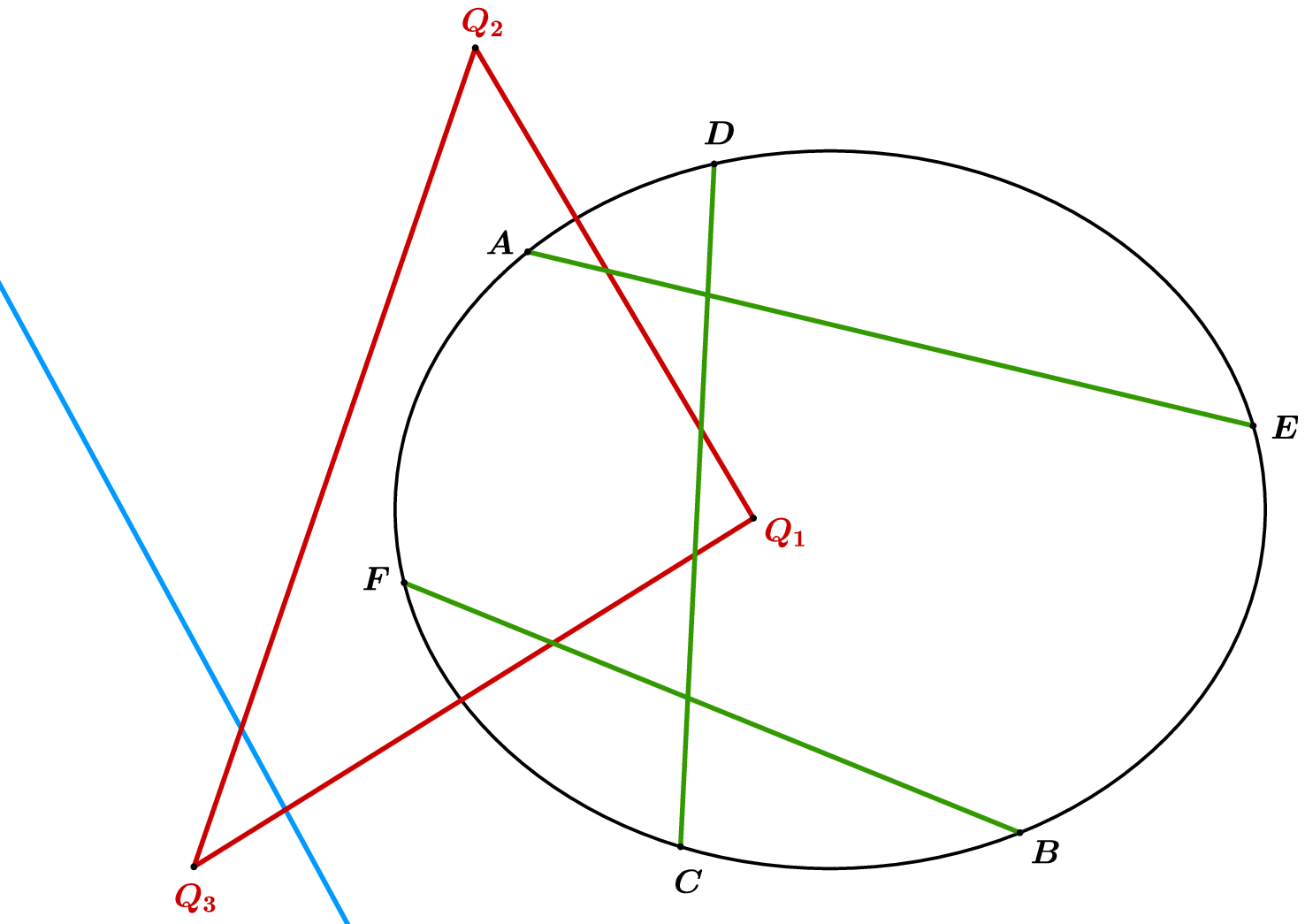}
\caption{\emph{Dramatis Personae} in the reconstruction} 
\label{diagram:fourpascals} 
\end{figure}

The relevant geometric elements are shown in
Diagram~\ref{diagram:fourpascals} on page~\pageref{diagram:fourpascals}. 
Since each Pascal corresponds to a $2 \times 3$ array (determined up to a
shuffling of rows and columns), its columns give a partition of the points $A, \dots, F$ 
into three sets of two elements each. For instance, any of the  
arrays in~(\ref{3pascals}) gives the partition  
\[ \{A, F\} \cup \{B, E\} \cup \{C, D\}. \] 

Now observe that the first three Pascals in~(\ref{list.fourpascals})
have been so chosen that they all lead to
the same partition, namely 
\begin{equation} 
\{A, E\} \cup \{C, D\} \cup \{B, F\}. 
\label{partition.first3arrays} \end{equation} 
This corresponds to the three green chords in
Diagram~\ref{diagram:fourpascals}. 
Let $Q_1$ denote the point $AB \cap EF$, which is common to $\ell_2$
and $\ell_3$. Similarly, let 
\begin{equation} Q_2 = \ell_3 \cap \ell_1 = AC \cap DE, \qquad 
Q_3 = \ell_1 \cap \ell_2 = BC \cap DF. 
\label{definitions.Qi} \end{equation}
Hence the line $Q_1Q_2$ is the same as $\ell_3$, and so on. Now, if we
switch the endpoints of all the three chords simultaneously; that is to
say, if we apply the transposition 
\[ (A \, E) \, (C \, D) \, (B \, F), \] 
then all the $Q_i$ remain unchanged and hence so do the first three
Pascals. In other words, each of the expressions $s_1, t_1, \dots,
s_3, t_3$ remains invariant if we make a simultaneous
substitution of variables $a \leftrightarrow e, c \leftrightarrow d, b
\leftrightarrow f$. It follows that no rational function of $s_1,
\dots, t_3$ can equal any of the variables $a, \dots, f$. 

The {\bf first stage} in the proof is to show that the next best outcome is achievable; that is
to say, the symmetric expressions 
\[ a+e, \quad a \, e, \quad b+f, \quad b \, f, \quad c+d, \quad c \, d \] 
are rational functions of $s_1, \dots, t_3$. In geometric terms (see
the top part of Diagram~\ref{diagram:fourpascals}), the red triangle $Q_1
Q_2 Q_3$ allows us to locate the three green chords, but 
we do not yet have sufficient information to label their
endpoints. The algebraic formulae which connect the red triangle to the green
chords are encoded in a transvectant identity. 

In the {\bf second stage}, we bring in the fourth Pascal $\ell_*$
(shown in blue) to break the symmetry. It is so chosen that each of
the three green chords passes through one of the cross-hair
intersections in $\ell_*$; for instance, $BF$ passes through the point 
$AD \cap BF$ on $\ell_*$. And now, another transvectant
identity allows us to get a linear equation for $a$ whose coefficients are rational functions in $s_1, \dots,
t_*$. This implies that $a$ itself is such a function, and a similar 
argument applies to $b, \dots, f$. This gives the required result. 

\subsection{Wernick's Problems} \label{section.Wernick} 
As an aside, we will point out the analogy between the
reconstruction problem for Pascals and Wernick's
problems~\cite{Schreck_etal, Wernick}. Given a triangle $ABC$ in the Euclidean
plane, one gets a large number of derived points such as the centroid,
the orthocentre or the three foots of perpendiculars. A
typical Wernick's problem asks whether the original triangle can be
reconstructed from a specific choice of three of the derived
points. Here are two sample results (see~\cite[p.~71]{Schreck_etal}): 
\begin{itemize} 
\item 
Given the centroid, orthocentre and the midpoint of any one side, the
original triangle is constructible. 
\item 
The original triangle is not constructible from the circumcentre,
orthocentre and the incentre. 
\end{itemize} 
It is clear that our main theorem is in this spirit, although the
specific geometric situation is different. The coordinates of
any of the derived points are often given by simple formulae in terms
of the coordinates of $A, B, C$. The analogous
formulae~(\ref{expressions.st}) in our case are more involved,
and hence the reconstruction is less straightforward. 

\section{Binary forms} \label{section.binary.forms} 
\subsection{} 
Let $\basefield$ denote the field $\QQ(a, b, c, d, e, f)$ of rational 
functions in the variables $a, \dots, f$. 
We will use $\basefield$ as our base field, so that  
any `scalar' will be assumed to belong to $\basefield$. Henceforth,
the projective plane $\P^2$ will be over $\basefield$. 

We will consider homogeneous forms in the variables $\ux = \{x_1,
x_2\}$. In a classical notation introduced by Cayley, 
$(z_0, z_1, \dots, z_n \cbrac x_1, x_2)^n$ 
stands for the degree $n$ form $\sum\limits_{i=0}^n \, z_i \,
\binom{n}{i} \, x_1^{n-i} \, x_2^i$. 

\subsection{Transvectants} 
Although the definition of a transvectant is \emph{prima facie}
technical, the concept arises naturally in invariant theory and representation
theory (see~\cite[Ch.~5]{Olver}).

Suppose that we are given two binary forms $G, H$ of
degrees $m,n$ respectively. For an integer $r \ge 0$, their $r$-th transvectant is defined to be 
\begin{equation} (G,H)_r = \frac{(m-r)! \, (n-r)!}{m! \, n!} \, 
\sum\limits_{i=0}^r \, (-1)^i \binom{r}{i} \, 
\frac{\partial^r G}{\partial x_1^{r-i} \, \partial x_2^i} \, 
\frac{\partial^r H}{\partial x_1^i \, \partial x_2^{r-i}} 
\label{trans.formula} \end{equation} 

This is a form of degree $m+n-2r$, unless it is identically zero.  
If 
\[ G = (g_0, g_1, g_2 \cbrac x_1, x_2)^2, \quad H = (h_0, h_1, h_2 \cbrac
x_1, x_2)^2, \] 
then it is easy to check that 
\[ 
\begin{aligned} 
(G, H)_1 & = (g_0 h_1 - g_1 h_0, \frac{1}{2} \, (g_0 \, h_2 - g_2 \,
h_0), g_1 \, h_2 - g_2 \, h_1 \cbrac x_1, x_2)^2, \; \text{and} \\ 
(G,H)_2 & = g_0 h_2 - 2 \, g_1 h_1 + g_2 h_0. 
\end{aligned} \] 
In general, the coefficients of $(G,H)_r$ are linear functions in the
coefficients of $G$ and $H$. The numerical factors in Cayley's
notation and~(\ref{trans.formula}) may seem unnecessary, but
experience has shown that they simplify the computations. 
\subsection{} 
Now the crucial step is to represent points and lines in
$\P^2$ by quadratic binary forms. (The reader may also refer  
to~\cite[\S 3]{Chipalkatti} where an identical set-up is used.) Let the nonzero 
quadratic form $G = (g_0, g_1, g_2 \cbrac x_1, x_2)^2$ represent the point $P_G = [g_0, g_1, g_2]$, as well as the line 
$L_G = \langle g_2, - 2 \, g_1, g_0 \rangle$. It is understood that
any nonzero scalar multiple of $G$ will represent the same point or
line. Now the following properties show that incidences
and joins are exactly mirrored by transvectants. 
\begin{Lemma} \rm With notation as above, 
\begin{enumerate} 
\item 
The point $P_G$ belongs to the line $L_H$, if and only if
$(G,H)_2=0$. 
\item 
The line joining the points $P_G$ and $P_H$ is $L_{(G,H)_1}$. 
\item 
The point of intersection of the lines $L_G$ and $L_H$ is $P_{(G,H)_1}$. 
\end{enumerate} \end{Lemma} 
All the proofs follow immediately from the definitions. The
point $P_G = [g_0, g_1, g_2]$ lies on $L_H = \langle
h_2, - 2 \, h_1, h_0 \rangle$ exactly when the dot product of the two
vectors is zero, which proves (1). The equation of the line joining
$P_G$ and $P_H$ is $\left| \begin{array}{ccc} z_0 & z_1 & z_2 \\ g_0 & g_1 & g_2 \\
             h_0 & h_1 & h_2 \end{array} \right|=0$, hence it is represented by $(G,H)_1$. The proof of (3) is
         similar. \qed 

\medskip 

The following result will be needed later. 
\begin{Lemma} \rm 
Two nonzero quadratic forms $G$ and $H$ are equal up to a scalar, if and only if
$(G,H)_1=0$. 
\label{lemma.jacobian} \end{Lemma} 
\proof    The forms are equal up to a scalar exactly when the 
matrix $\left[ \begin{array}{ccc} g_0 & g_1 & g_2 \\ h_0 & h_1 & h_2 \end{array}
\right]$ has rank one, i.e., exactly when all of its minors are
zero. This is equivalent to the vanishing of all the coefficients of
$(G,H)_1$. \qed 

\medskip 

The advantage of using transvectants is that there are well-developed
tools for manipulating them, namely, a symbolic calculus (see~\cite{GraceYoung, Olver})
as well as a graphical calculus (see~\cite[\S2]{Abdesselam}).
This is especially useful when one encounters transvectants whose
components are themselves transvectants. 

\subsection{} 
The conic $\conic$ consists of those points $P_G$ such that
\[ (G,G)_2= 2 \, (g_1^2 - g_0 \, g_2)=0. \] 
These are the nonzero forms $G$ which can be written as 
squares of linear forms up to a scalar. Define six linear forms 
\[ a_\ux = x_1 + a \, x_2, \quad b_\ux = x_1 + b \, x_2, \quad \dots
\quad f_\ux = x_1 + f \, x_2,  \] 
and fix the points $A = P_{a_\ux^2}, \dots, F = P_{f_\ux^2}$ on $\conic$. Let 
\begin{equation} \lambda_i = (t_i, - \frac{s_i}{2}, 1 \cbrac x_1, x_2)^2, \qquad i
=1,2,3,* \label{definition.lambda_i} \end{equation} 
denote the quadratic forms which represent the Pascals $\ell_i$. All
of this agrees with the notational conventions in
section~\ref{section.maintheorem.refined}. 

\begin{figure}
\includegraphics[width=8cm]{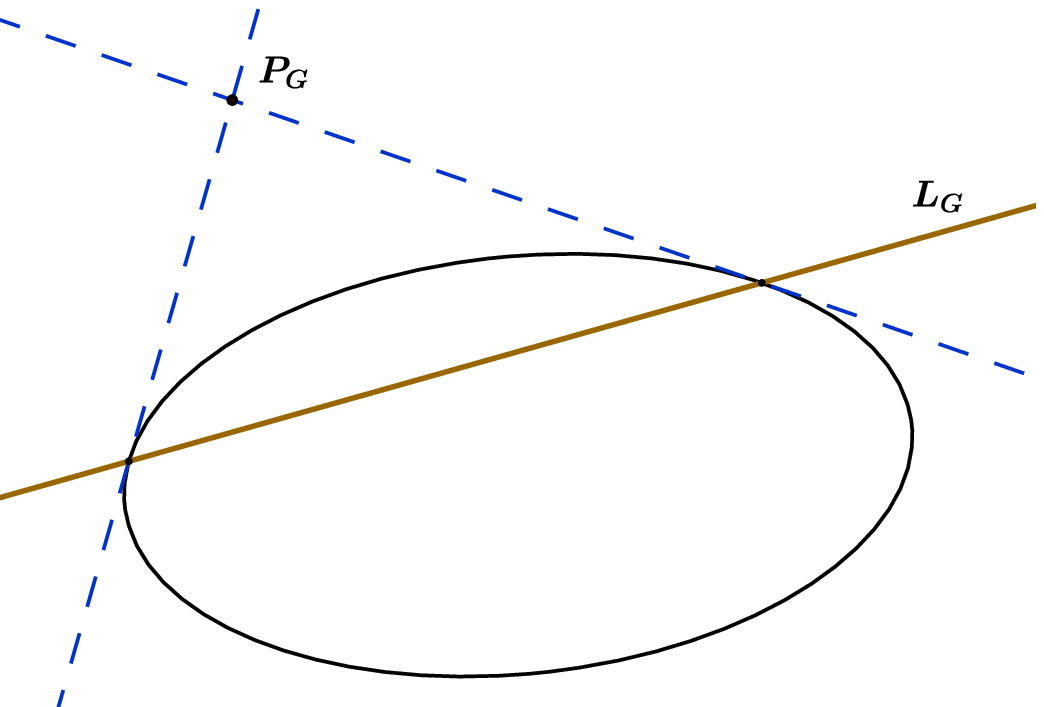}
\includegraphics[width=7cm]{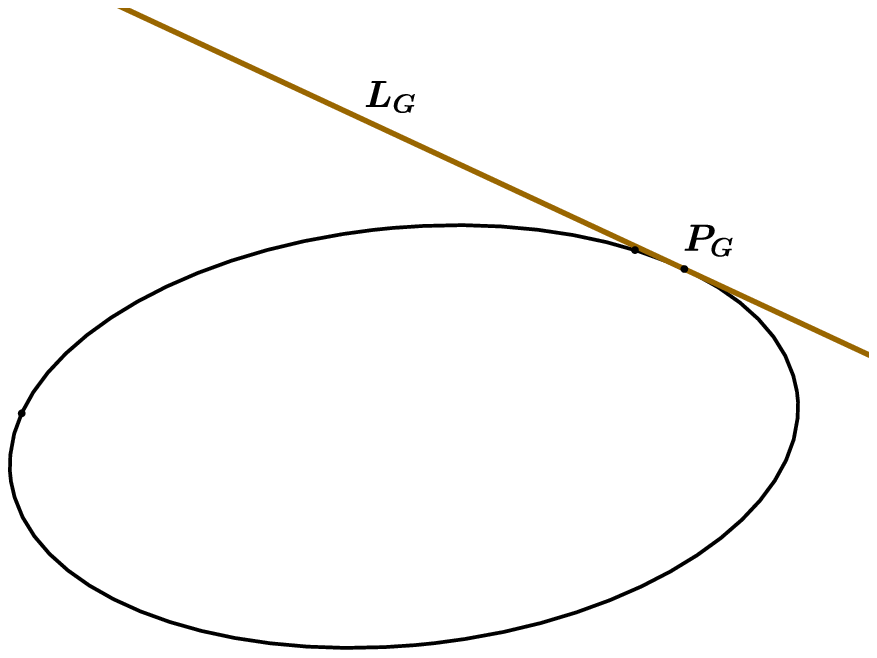}
\caption{The pole-polar relation} \label{diagram:polepolar}
\end{figure}

The following lemma is helpful in completing the geometric picture, but it will
not be needed elsewhere (see Diagram~\ref{diagram:polepolar}). 

\begin{Lemma} \rm 
Let $G$ denote a nonzero quadratic form. Then $L_G$ is the
polar line of $P_G$ with respect to $\conic$. In particular, 
\[ \text{$P_G$ lies on $L_G$} \iff \text{$P_G$ lies on $\conic$} \iff 
\text{$L_G$ is tangent to $\conic$}. \] 
\end{Lemma} 
The proof is left to the reader. \qed

\subsection{} 
For instance, the line $AB$ is represented by the form $(a_\ux^2,
b_\ux^2)_1 = (b-a) \, a_\ux \, b_\ux$, or after ignoring the scalar, just by $a_\ux \, b_\ux$. 
It follows that the points $Q_1, Q_2, Q_3$ in
section~\ref{section.proof.overview} are respectively represented by
the quadratic forms 
\begin{equation} \pi_1 = (a_\ux b_\ux, e_\ux f_\ux)_1, \quad 
\pi_2 = (a_\ux c_\ux, d_\ux e_\ux)_1, \quad 
\pi_3 = (b_\ux c_\ux, d_\ux f_\ux)_1. 
\label{formulae.pi_j} \end{equation} 
Since $Q_1 = \ell_2 \cap \ell_3$ etc, they are also respectively
represented by
\begin{equation} \mu_1 = (\lambda_2, \lambda_3)_1, \quad 
\mu_2 = (\lambda_3, \lambda_1)_1, \quad 
\mu_3 = (\lambda_1, \lambda_2)_1. 
\label{formulae.mu_j} \end{equation} 
This implies that $\mu_i$ and $\pi_i$ are equal up to a multiplicative
scalar in $\basefield$. It is clear that the coefficients of $\mu_i$
are rational functions in $s_1, \dots, t_3$. 

\section{The proof of the main theorem} \label{section.main.proof} 
\subsection{The first stage} \label{section.firststage}
For any quadratic forms $U, V, W$, define 
\[ \psi(U, V, W) = 6 \, (U, V W)_2 - U  (V, W)_2,  \] 
which is also a quadratic form. 

\begin{Proposition} \rm 
We have an identity 
\begin{equation} \psi(\pi_3, \pi_1, \pi_2) = \Phi \times a_\ux \,
  e_\ux,  
\label{ae.identity} \end{equation} 
where $\Phi$ is a polynomial in $a, \dots, f$. 
\label{proposition.ae} \end{Proposition} 

The proof will be given in section~\ref{section.graphical.proofs}
using the graphical calculus, 
but the rationale behind the proposition can be explained without
it. The right-hand side of~(\ref{ae.identity}) represents the line $AE$. Since $\mu_i$ is proportional to
$\pi_i$, the left-hand side is proportional to $\psi(\mu_3,\mu_1,\mu_2)$. Hence
the identity implies that $AE$ can be represented by a form 
\[ \lambda_{AE} = (\alpha_{AE}, \beta_{AE}, 1 \cbrac x_1, x_2)^2, \] 
where $\alpha_{AE}, \beta_{AE}$ are rational functions of $s_1, \dots,
t_3$. We can similarly write down $\lambda_{CD}$ and $\lambda_{BF}$
representing the other two green chords in Diagram~\ref{diagram:fourpascals}. 
The exact expression for $\Phi$ will be found in the course of proving the identity, but it is immaterial
to the main theorem. 

Formula~(\ref{ae.identity}) was initially obtained by some 
calculated guesswork guided by intuition. Since the construction of
Pascals is synthetic, if it is at all 
possible to pass from the red triangle to the green chords, then the connecting formula can be
plausibly written in terms of transvectants. Since the letters $a,e$ enter symmetrically into the expressions
for $\pi_1, \pi_2$, the formula should respect this structure as
well. Now the correct definition of $\psi$ is determined by a graphical
calculation, in which the initial intuition is buttressed by a formal proof. 

A direct calculation shows that 
\[ \alpha_{AE} = \frac{ s_1^2 \, s_3 \, t_2 \, t_3 -s_1^2 \, s_2 \,
  t_3^2 + \text{$16$ similar terms}}{s_1^2 \, s_2 \, t_2 + s_1 \, s_2
  \, s_3 \, t_1 + \text{$16$ similar terms}}, \] 
with a similar expression for $\beta_{AE}$. Thus~(\ref{ae.identity})
serves as a compact shorthand for a lengthy and complicated formula. 

\subsection{The second stage} 
We now use the fourth Pascal $\ell_*$. 
Recall that a point on $\conic$ is represented by the square of a
linear form which is well-defined up to a scalar. Thus $A$ comes from 
$a_\ux$, where we are hoping to solve for $a$ in terms of $s_1, t_1, \dots, s_*, t_*$. 
There are two ways of expressing $D$ in terms of $A$, and their comparison
will lead to a set of equations for $a$.

\begin{figure}
\includegraphics[width=10cm]{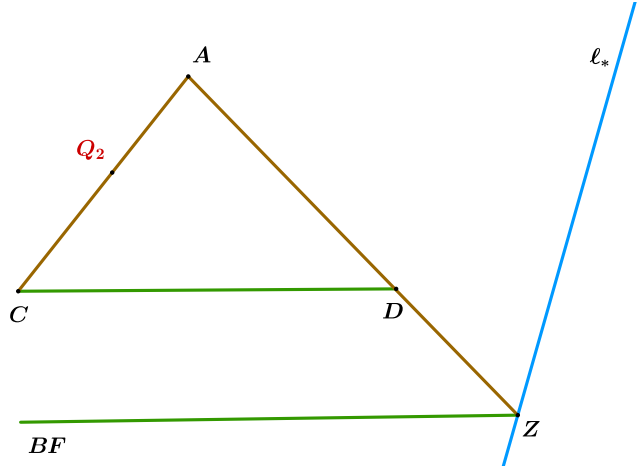}
\caption{The double dependence of $D$ on $A$ } 
\label{diagram:pointZ} 
\end{figure}  

Diagram~\ref{diagram:pointZ} shows the geometric elements needed
in the second step. 
\begin{enumerate} 
\item 
Since the point $Q_2$ is on $AC$, the line $AQ_2$ is the same as
$AC$. Now $AQ_2$ is represented by 
\[ (\mu_2, a_\ux^2)_1 = a_\ux \, (\mu_2, a_\ux)_1. \] 
Hence $C$ comes from the linear form $(\mu_2, a_\ux)_1$, and thus $D$
comes from 
\begin{equation}  \frac{\lambda_{CD}}{(\mu_2,a_\ux)_1}. 
\label{form.D1} \end{equation} 
\item 
The Pascal $\ell_*$ passes through $Z=AD \cap BF$, which implies that 
$Z$ is represented by $(\ell_\star, \lambda_{BF})_1$. Hence $AZ$, which is the same as $AD$, is represented by 
\[ ((\ell_*, \lambda_{BF})_1, a_\ux^2)_1 = a_\ux \, ((\ell_*,
\lambda_{BF})_1, a_\ux)_1. \] 
Thus $D$ also comes from 
\begin{equation} 
((\ell_*, \lambda_{BF})_1, a_\ux)_1. \label{form.D2} \end{equation}
\end{enumerate} 
The two linear forms in (\ref{form.D1}) and (\ref{form.D2}) must coincide up to a scalar. This gives the identity 
\[ \lambda_{CD} = \text{scalar} \times (\mu_2, a_\ux)_1 \, \times  
((\ell_*, \lambda_{BF})_1, a_\ux)_1. \] 
If we write 
\[ U = \lambda_{CD}, \quad V = \mu_2, \quad W =
(\ell_*,\lambda_{BF})_1, \] 
then, by Lemma~\ref{lemma.jacobian}, this is equivalent to 
\[ (U, (V,a_\ux)_1 \, (W, a_\ux)_1)_1 =0. \] 

The following transvectant identity allows us to rewrite this in such a way that
we can extract a set of equations for $a$. 
\begin{Proposition} \rm 
For arbitrary quadratic forms $U, V, W$ and linear form $a_\ux$, we
have an identity 
\begin{equation} 
(U, (V,a_\ux)_1 \, (W, a_\ux)_1)_1 = (M, a_{\ux}^2)_2 + (N, a_{\ux}^2)_1,  
\label{ax.identity} \end{equation} 
where 
\[ M = \frac{1}{2} \, (U, W)_1 \, V + \frac{1}{2} \, (U, V)_1 \, W, 
\qquad 
N = - \frac{1}{2} \, (U, V W)_2 - \frac{1}{6} \, U  (V,
W)_2. \] 
\label{proposition.ax} \end{Proposition} 

The proof will be given in
section~\ref{section.graphical.proofs}. The purpose of the identity 
is to `package' the known quantities $U, V, W$ into $M$ and $N$, so as to separate them from the unknown
quantity $a$.

\subsection{} 
Now write 
\[ M = (m_0, m_1, m_2,m_3,m_4 \cbrac x_1,x_2)^4, \quad \text{and}
\quad N = (n_0, n_1, n_2 \cbrac (x_1,x_2)^2. \] 
The coefficients of $U, V, W$ are rational functions of $s_1,
\dots, t_*$, hence so are all the $m_i$ and $n_i$. 
The right-hand side of (\ref{ax.identity}) can be expanded as 
$(r_0, r_1, r_2 \cbrac x_1, x_2)^2$, where each $r_i$ is quadratic in $a$. Since this must vanish identically,
we get three quadratic equations $r_0=r_1=r_2=0$ for $a$. 
A straightforward expansion shows that they can be written as 
\[ \left[ \begin{array}{ccc} 
m_2 - n_1 & n_0 - 2 \, m_1 & m_0 \\ 
2 \, m_3 - n_2 & - 4 \, m_2 & 2 \, m_1 +  n_0 \\ 
m_4 & - 2 \, m_3 - n_2 & m_2 + n_1 
\end{array} \right] 
\left[ \begin{array}{c} 1 \\ a \\ a^2 
\end{array} \right] = 0. \] 
Let $Z = (z_{ij})$ denote the $3 \times 3$ matrix on the left; e.g., $z_{12} = n_0 - 2m_1$. 
Now, for instance, we can use its first two rows to solve for 
$a$, which gives 
\[ a = - \frac{\left| \begin{array}{cc} z_{11} & z_{13} \\ z_{21} &
                                                                    z_{23} \end{array} \right|}{
\left| \begin{array}{cc} z_{12} & z_{13} \\ z_{22} &
                                                     z_{23} \end{array}
                                                 \right|}. \] 
This proves that $a$ is a rational function of $s_1, \dots,
t_*$. Since $e_\ux$ is a constant multiple of
$\frac{\lambda_{AE}}{a_\ux}$, the same follows for $e$. The Pascal
$\ell_*$ passes through the points $CD \cap BE, AE \cap CF$ which
respectively lie on the green chords $CD, AE$. Hence the same argument
as in the second stage gives the result for $b, c, d, f$. This proves the main
theorem, assuming Propositions~\ref{proposition.ae} and~\ref{proposition.ax}. \qed

The passage 
\[ \{s_1, t_1, \dots, s_*, t_*\} \Rightarrow \{a, \dots, f\} \] 
goes through two complicated algebraic identities neither of which has any 
obvious geometric content. Thus our reconstruction is not
`synthetic' in the classical sense of the word. We do not know of any 
natural ruler-and-compass type construction which begins with the Pascals 
and ends with the sextuple. It would be 
interesting to find one. 

\subsection{} 
The main theorem is valid over any field of characteristic zero, since
the choice of $\QQ$ plays no essential role in the proof. Moreover, 
the only numerical coefficients which appear in the proof are 
$2, 4, 6$ and $\frac{3}{4}$. All of these are \emph{defined and
  nonzero} as long as the base field has characteristic $\neq 2,
3$, and hence the theorem remains valid over such a field. It would be
interesting to have a similar theorem when the characteristic is
either $2$ or $3$. 

\subsection{} We have programmed the entire procedure in {\sc Maple}
in order to ensure against the possibility of error. For instance,
suppose that 
\[ a =7, \quad b = -3, \quad c = 2, \quad d = 5, \quad e = -4, \quad f
=1. \] 
Then the Pascals are 
\[ \begin{array}{lll} 
\lambda_1 = (\frac{5}{36}, \frac{37}{72},1 \cbrac x_1,x_2)^2, & & 
\lambda_2 = (-\frac{49}{349}, \frac{42}{349},1 \cbrac x_1,x_2)^2, \\ 
\lambda_3 = (-\frac{1}{16}, -\frac{33}{544},1 \cbrac x_1,x_2)^2, & & 
\lambda_* = (\frac{7}{74}, \frac{21}{148},1 \cbrac x_1,x_2)^2. 
\end{array} \] 
Now if we follow the recipe given above, the result is 
\[ a = 
\frac{5^7.7^2.11^{12}.13^4.29^{14}}{5^7.7.11^{12}.13^4.29^{14}} =
7 \] 
as expected, and similarly for the remaining variables. We have
done a similar verification on several such examples. 

\subsection{}
The theme of this paper is related to the Galois (or monodromy) group of Pascal 
  lines in the sense of~\cite{Harris}. We explain this in 
brief. 

Assume the base field to be $\complex$. Write 
\[ (T-a) \, (T-b) \, \dots (T-f) = T^6 - s_1 \, T^5 + s_2 \, T^4 - s_3 
\, T^3 + s_4 \, T^2 - s_5 \, T +s_6, \] 
where $s_1, \dots, s_6$ are the elementary symmetric functions in $a,
\dots, f$. Let $\unord$ denote the space of \emph{unordered} six 
points on a conic. In fact $\unord$ is birational to $\text{Sym}^6 
\conic \simeq \P^6$, and its field of rational functions may be identified with 
\[ \FF = \complex(s_1, \dots, s_6). \] 
We have a $60$-$1$ cover $\cover \lra \unord$, where the fibre over an 
unordered sextuple corresponds to its collection of $60$ Pascals. If 
$\langle 1, s_i, t_i \rangle, 1 \le i \le 60$ are the line coordinates 
of the Pascals, then the field of rational functions of $\cover$ is 
$\FF(s_1, t_1, \dots, s_{60}, t_{60})$. However, the inclusion 
\begin{equation} \FF(s_1, t_1, \dots, s_{60}, t_{60}) \subseteq 
  \complex(a, \dots, f) 
\label{field.inclusion} \end{equation} 
is actually an equality by our main theorem. Hence we have the 
following: 
\begin{Proposition} \rm 
The Galois group 
\[ \Gal(\cover/\unord) \simeq \Gal(\complex(a,\dots,f)/\complex(s_1, \dots,
s_6)) \] 
is isomorphic to the symmetric group on six letters. 
\end{Proposition} 

It should be clarified that this result cannot be considered original 
to this paper. The fact that~(\ref{field.inclusion}) is an equality is 
already implicit in Pedoe's proof in~\cite{Pedoe}, although it is not 
so stated there.

\subsection{Optimal subsets} 
Let $X$ be an arbitrary $n$-element subset of the sixty Pascals, with line 
coordinates 
\[ \langle 1, s^{(i)}, t^{(i)} \rangle, \quad 1 \le i \le n. \] This gives an inclusion of fields 
\begin{equation} 
\underbrace{\QQ(s^{(1)}, t^{(1)}, \dots, s^{(n)},
  t^{(n)})}_{\basefield_X}   \subseteq \basefield. \label{field.extension} \end{equation}
Let us say that the set $X$ is \emph{adequate} if equality holds; 
this is equivalent to saying that each variable is a rational 
function in $s^{(1)}, \dots, t^{(n)}$. Furthermore, let us say that 
$X$ is \emph{optimal} if it is adequate and no proper subset of $X$ is 
adequate. 

\begin{Proposition} \rm 
The set of Pascals given in the main theorem is optimal. 
\end{Proposition} 

\proof It is clear that $\{ \ell_1, \ell_2, \ell_3\}$ is not adequate,
in fact~(\ref{field.extension}) is a quadratic extension in this
case. If we take $X = \{\ell_1, \ell_2, \ell_*\}$, then a Maple
computation shows that~(\ref{field.extension}) is a degree $12$
extension, and hence $X$ is not adequate. (It would be better to
have a more conceptual and less computational proof, but we cannot
find one.) 

Now observe that the permutation $(A \, F) \, (B \, E) \, (C \, D)$
leaves $\ell_1, \ell_*$ unchanged, and interchanges $\ell_2$ and
$\ell_3$. Hence the same result follows for $\{\ell_1, \ell_3,
\ell_*\}$. Finally, the permutation $(A \, D) \, (B \, F) \, (C \, E)$
interchanges $\ell_1, \ell_3$ and leaves $\ell_2, \ell_*$ unchanged,
which proves that $\{\ell_2, \ell_3, \ell_*\}$ is not adequate. This
completes the proof. \qed 

Since $\basefield$ has transcendence degree $6$ over $\QQ$, any
adequate subset must have at least $3$ elements. It would be of interest
to know whether there exists an adequate $3$-element subset, which must
then be necessarily optimal. We have not succeeded in finding any. 

On the other hand, given an arbitrary subset of (three or more) Pascals,
it is not at all obvious how to decide whether it is adequate. Thus
there is a large number of Wernick-Pascal type reconstruction problems which remain
open. It is a matter of speculation whether transvectant identities of
some sort will play a role in their solution. 

\subsection{} There are geometric obstructions which prevent certain
sets from being adequate. Consider the set $X$ consisting of Pascals 
\[ \pascal{A}{B}{C}{F}{E}{D},\quad 
\pascal{A}{B}{C}{D}{F}{E},\quad 
\pascal{A}{B}{C}{E}{D}{F}, \] 
where the top row is held constant and the bottom row undergoes a
cyclic shift. Steiner's theorem says that
these three Pascals are concurrent. If $\langle 1, s^{(i)}, t^{(i)} \rangle,
i=1,2,3$ denote their line coordinates, then the determinant 
$\left| \begin{array}{ccc} 1 & s^{(1)} & t^{(1)} \\ 1 & s^{(2)} & t^{(2)} 
\\ 1 & s^{(3)} & t^{(3)} \end{array} \right| =0$. Hence $\basefield_X$ has
transcendence degree at most $5$ over\footnote{It can be shown to be exactly
  $5$, but this is not needed for the conclusion.} $\QQ$, and $X$ cannot be
adequate. Rather similarly, Kirkman's theorem says that the Pascals 
\[ \pascal{A}{B}{C}{F}{E}{D},\quad 
\pascal{A}{D}{F}{C}{E}{B}, \quad 
\pascal{A}{C}{F}{E}{B}{D}, \] 
are concurrent, and then the same conclusion follows. The reader will find 
a proof of either theorem in Salmon's notes referred to above. 
                                                                 
\section{Transvectant identities} 
\label{section.graphical.proofs} 
In this section we will prove Propositions~\ref{proposition.ae}
and~\ref{proposition.ax}. The proofs rely upon the graphical
formalism\footnote{It has a close affinity to the classical
  symbolic calculus as practiced by the German school of invariant
  theorists in the nineteenth century (cf.~\cite{Clebsch,
    GraceYoung,KungR}). The bibliography of~\cite{Abdesselam} contains 
  several more references to this circle of ideas.} 
developed in~\cite[\S2]{Abdesselam}. 

\subsection{} 
We will first rewrite Proposition~\ref{proposition.ae}
in more general and precise form. Consider six general linear forms 
$a_\ux=a_1x_1+a_2x_2$, $b_\ux=b_1x_1+b_2x_2$, \ldots,
$f_\ux=f_1x_1+f_2x_2$, where a letter such as `$a$' stands for a pair of variables $(a_1,a_2)$ instead of a single one.
We will also use the classical bracket notation $(ab)=a_1b_2-a_2 b_1$
for $2\times 2$ determinants, and similarly for $(cd)$, $(bf)$, etc.

Write 
\[ U=(b_\ux c_\ux,d_\ux f_\ux)_1, \quad 
V=(a_\ux c_\ux,d_\ux e_\ux)_1, \quad W=(a_\ux b_\ux,e_\ux f_\ux)_1, \] 
and~$\psi(U,V,W)=6(U,V W)_2-U(V,W)_2$. 
Define 
\begin{equation} \exprS=(da)(fc)(eb)-(ce)(bd)(af). 
\label{defn.exprS} \end{equation} 
\begin{Proposition} \rm 
With notation as above, we have 
\[ \psi(U,V,W)=\Phi\ a_\ux e_\ux, \]
where 
$ \Phi=\frac{3}{4}(cd)(bf)\ \mathcal{S}$. 
\label{prop31} \end{Proposition}

\begin{Remark}\label{sym.remark} \rm 
The expression $\exprS$ has the following invariance property.  Let
$J$ denote the operation of making a simultaneous exchange of letters
$a \leftrightarrow b, e \leftrightarrow f$. Now $\exprS$ remains
invariant under the action of $J$, since the bracket factors $(da), (fc), (eb)$ are
respectively taken to $(db), (ec), (fa)$ and conversely. 
Similarly, let $K$ and $L$ respectively denote the operations 
\[ b \leftrightarrow c, f \leftrightarrow d, \quad \text{and} \quad 
a \leftrightarrow e, b\leftrightarrow f, c \leftrightarrow d. \] 
Then $K$ also leaves $\exprS$ invariant, whereas $L$ changes it to
$-\exprS$. The subgroup generated by $J,K$ and $L$ inside the permutation group on letters $a,
\dots, f$, is isomorphic to $\SG_3 \times \ZZ_2$. 
\end{Remark}

\begin{Lemma}\label{psilem} \rm 
We have the more symmetric rewriting
\[
\psi(U,V,W)=3\left[
(U,V)_2 W+(U,W)_2 V-(V,W)_2 U
\right]\ .
\]
\end{Lemma}

\proof
Using the graphical formalism of~\cite[\S2]{Abdesselam}, we can write 
\begin{equation}
(U,VW)_2=\ \ 
\parbox{2.2cm}{
\psfrag{U}{$\scriptstyle{U}$}
\psfrag{V}{$\scriptstyle{V}$}
\psfrag{W}{$\scriptstyle{W}$}
\psfrag{x}{$\scriptstyle{x}$}
\includegraphics[width=2.2cm]{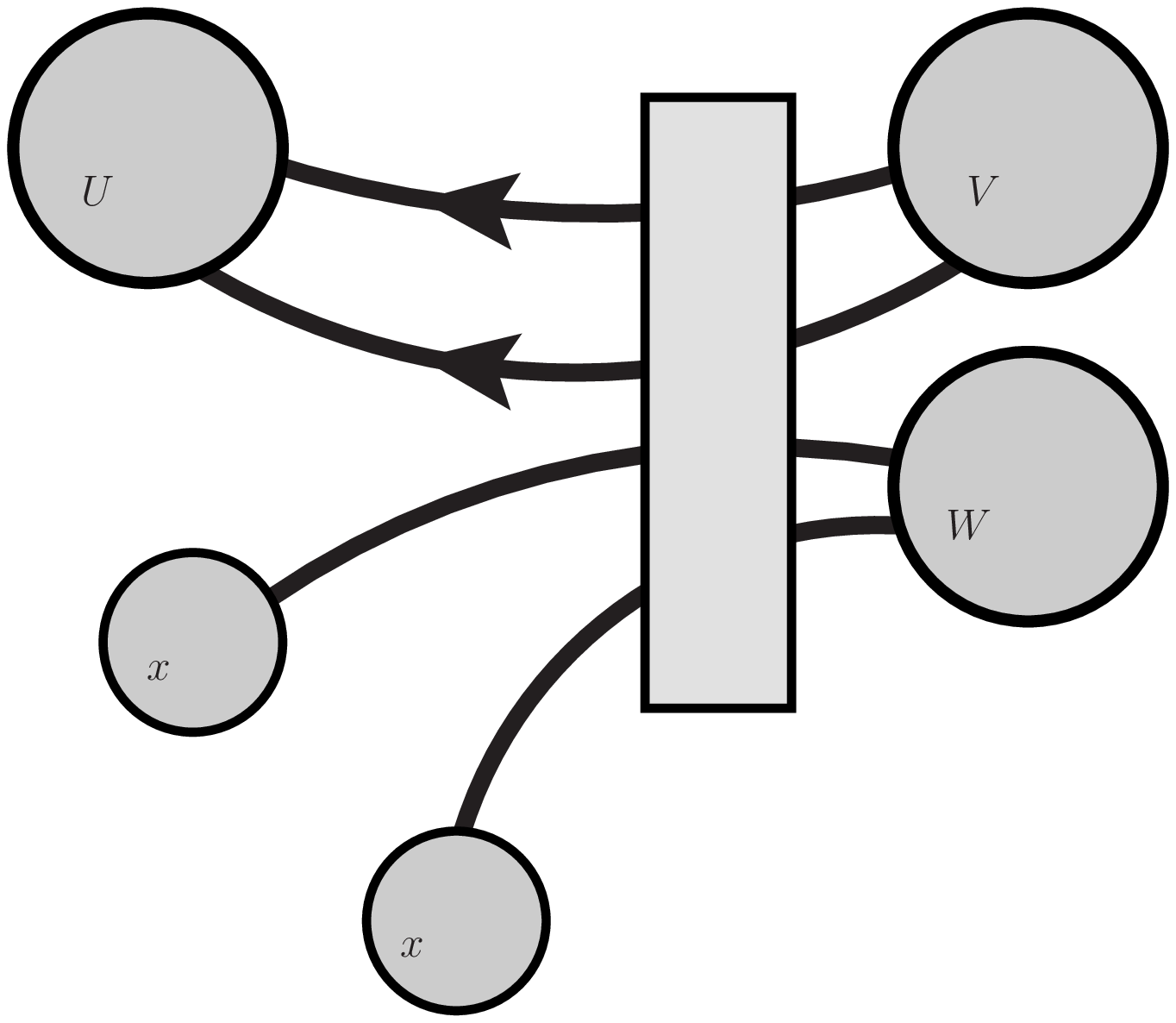}}
\ \ =\frac{1}{4!}\left[
4(U,V)_2 W+4(U,W)_2 V+16\ \ 
\parbox{2.5cm}{
\psfrag{U}{$\scriptstyle{U}$}
\psfrag{V}{$\scriptstyle{V}$}
\psfrag{W}{$\scriptstyle{W}$}
\psfrag{x}{$\scriptstyle{x}$}
\includegraphics[width=2.5cm]{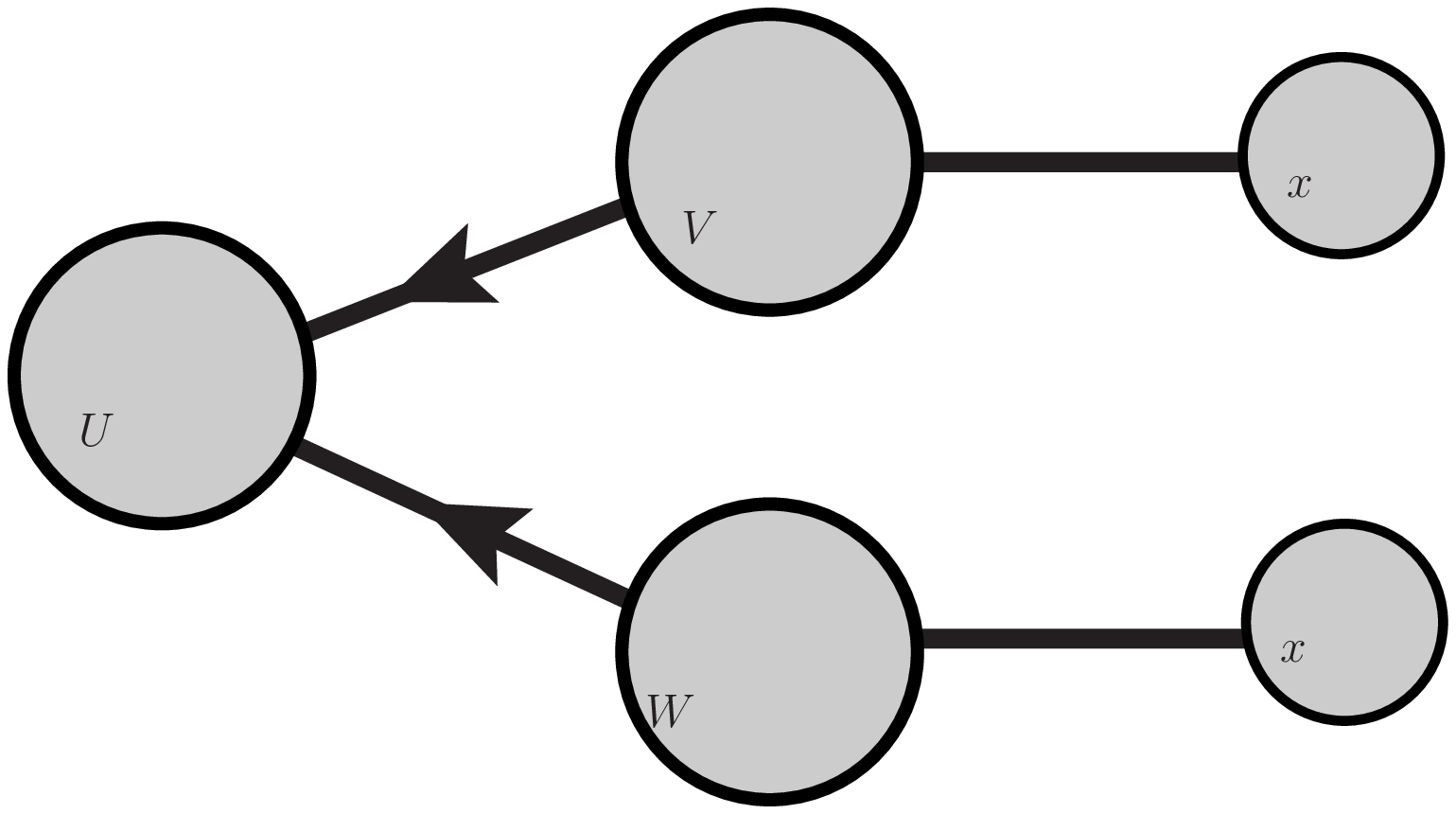}}
\ \ \right]
\label{psiexpeq}
\end{equation}
by expanding the normalized $\SG_4$ symmetrizer (represented by the grey rectangle).
We will use the notation
\[
\{V\rightarrow U\leftarrow W\}=\ \parbox{2.5cm}{
\psfrag{U}{$\scriptstyle{U}$}
\psfrag{V}{$\scriptstyle{V}$}
\psfrag{W}{$\scriptstyle{W}$}
\psfrag{x}{$\scriptstyle{x}$}
\includegraphics[width=2.5cm]{pic2.eps}}\ \ .
\]
Inserting the matrix identity $\epsilon\epsilon^{\rm T}=I$ (where $\epsilon$ is the $2\times 2$ antisymmetric matrix
with $\epsilon_{12}=1$ represented by the arrows), 
and using the Grassmann-Pl\"{u}cker (GP) relation where indicated by
the dotted line, we have
\[
\{V\rightarrow U\leftarrow W\}=\ 
\parbox{2.5cm}{
\psfrag{U}{$\scriptstyle{U}$}
\psfrag{V}{$\scriptstyle{V}$}
\psfrag{W}{$\scriptstyle{W}$}
\psfrag{x}{$\scriptstyle{x}$}
\psfrag{G}{$\scriptstyle{GP}$}
\includegraphics[width=2.5cm]{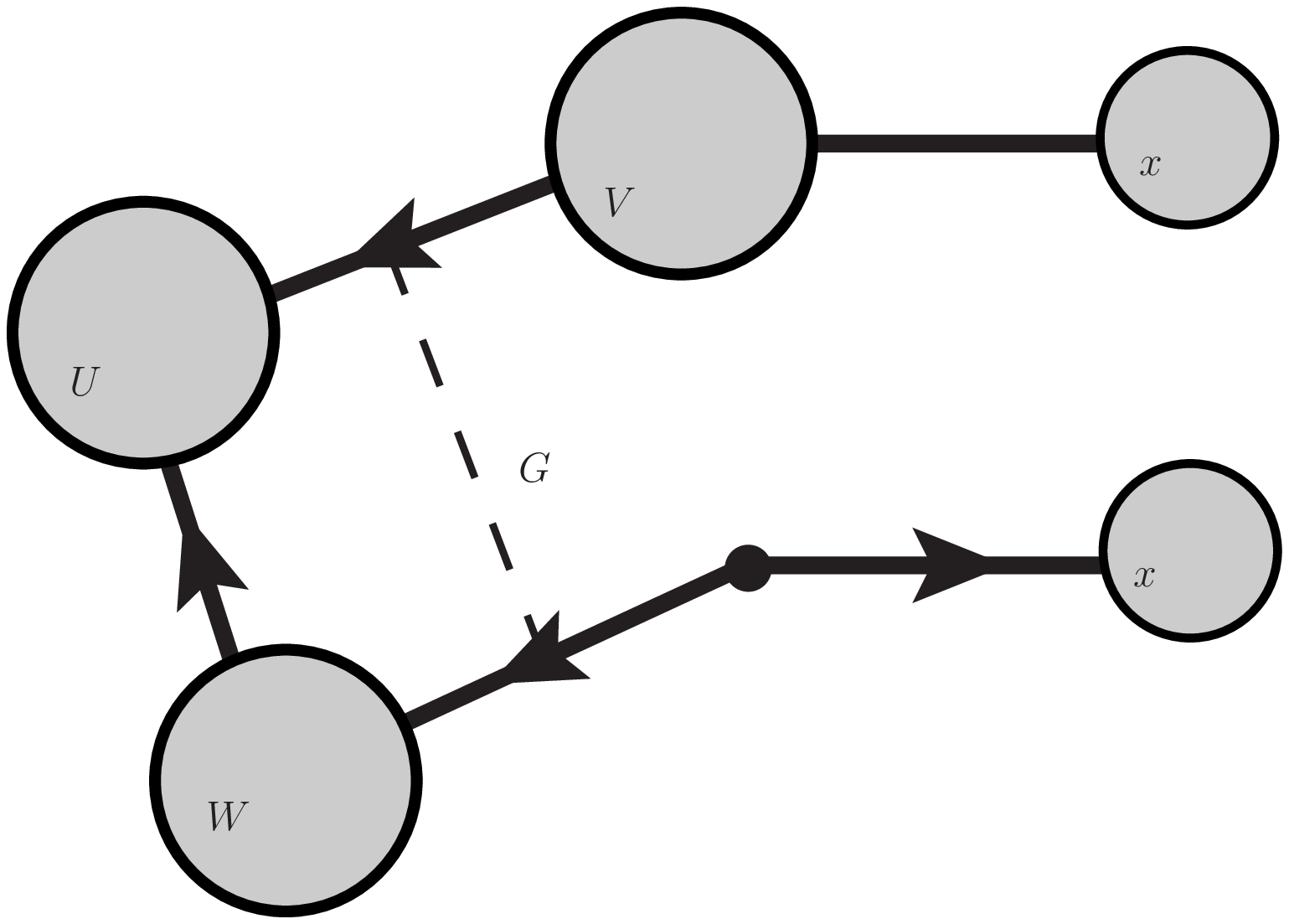}}
\ =\ 
\parbox{2.5cm}{
\psfrag{U}{$\scriptstyle{U}$}
\psfrag{V}{$\scriptstyle{V}$}
\psfrag{W}{$\scriptstyle{W}$}
\psfrag{x}{$\scriptstyle{x}$}
\includegraphics[width=2.5cm]{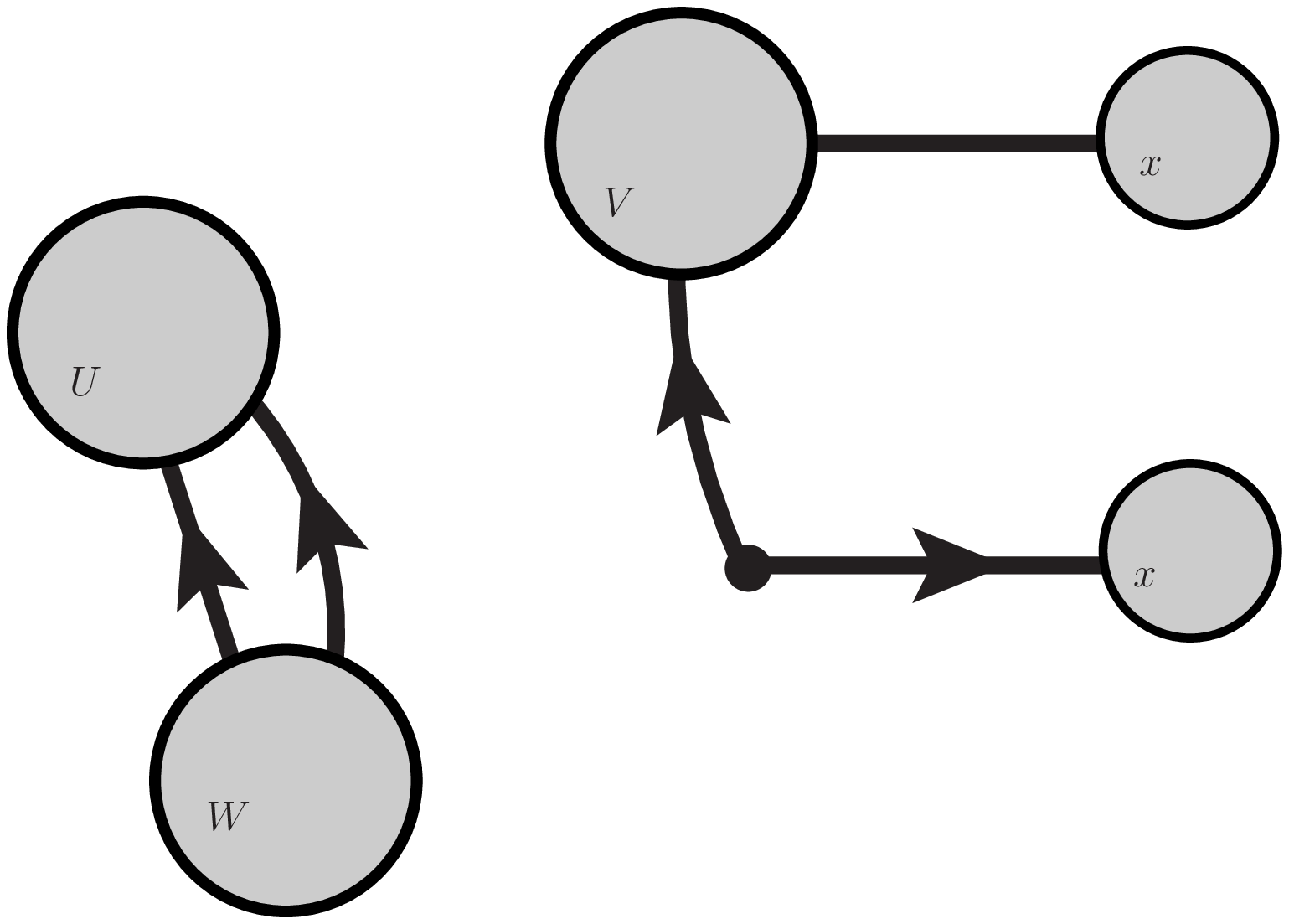}}
\ -\ 
\parbox{2.5cm}{
\psfrag{U}{$\scriptstyle{U}$}
\psfrag{V}{$\scriptstyle{V}$}
\psfrag{W}{$\scriptstyle{W}$}
\psfrag{x}{$\scriptstyle{x}$}
\includegraphics[width=2.5cm]{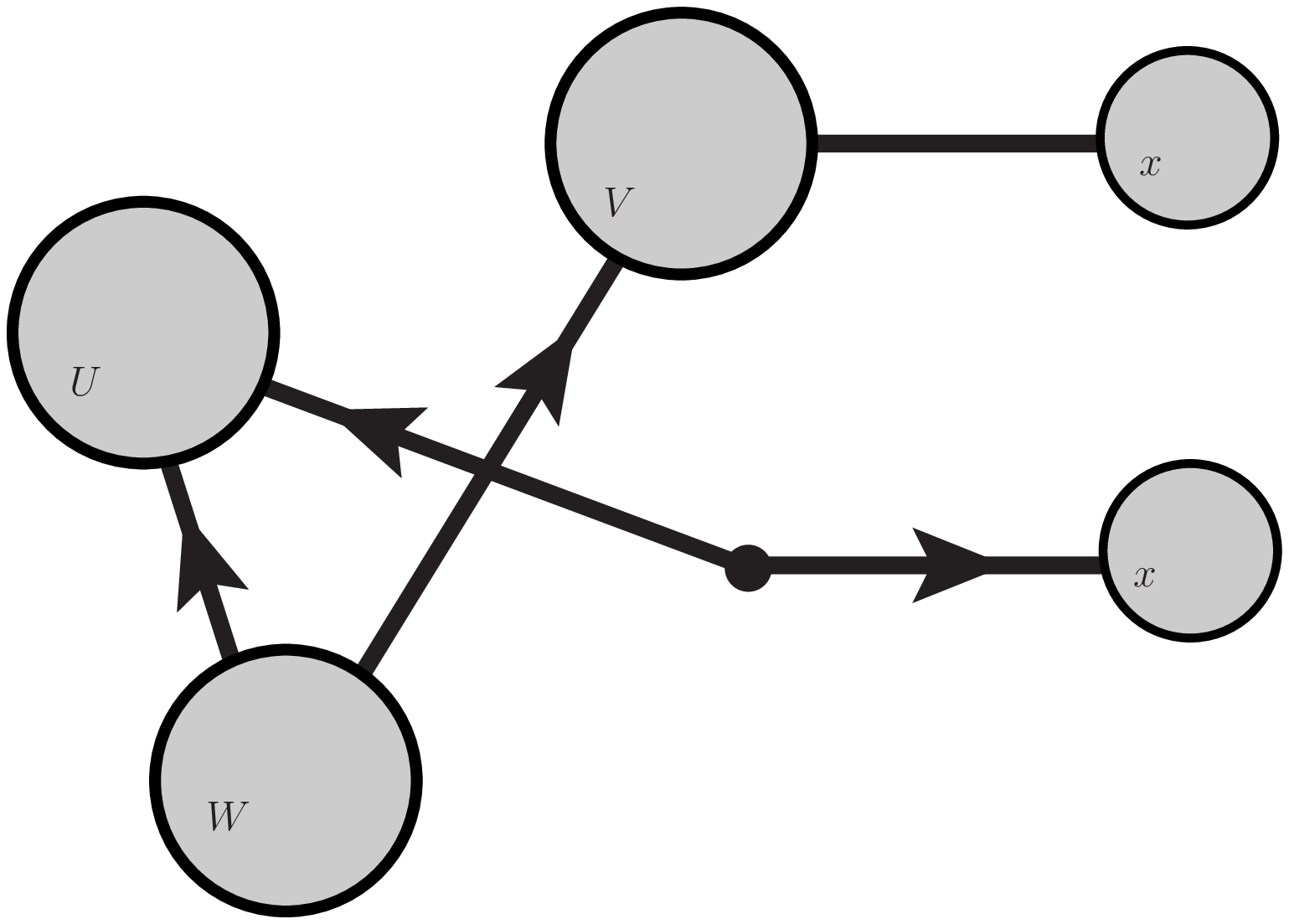}}\ \ ,
\]
i.e.,
\begin{equation}
\{V\rightarrow U\leftarrow W\}+\{U\rightarrow W\leftarrow V\}=(U,W)_2V\ .
\label{twochaineq}
\end{equation}
Permuting $U$, $V$ and $W$ in the last identity
gives three equations. They can be written in matrix form as
\[ \left(\begin{array}{c}
(V,W)_2 U\\ (U,W)_2 V\\ (U,V)_2 W \end{array} \right) 
=\left(\begin{array}{ccc}
0 & 1 & 1\\
1 & 0 & 1\\
1 & 1 & 0
\end{array}\right)
\ \left(
\begin{array}{c}
  \{V \rightarrow U \leftarrow W]\}  \\
 \{U\rightarrow V\leftarrow W\} \\
 \{U\rightarrow W\leftarrow V\}
\end{array}
\right)\ .
\]
By inverting this matrix, we get
\begin{equation}
\{V\rightarrow U\leftarrow W\}=\frac{1}{2}\left[
-(V,W)_2 U+(U,W)_2 V+(U,V)_2 W
\right]\ .
\label{onechaineq}
\end{equation}
After substituting back in (\ref{psiexpeq}) and simplifying, we get the required expression. \qed

The next lemma will be useful in the calculation of $\psi$. 
\begin{Lemma}\label{fourtotwo} \rm 
We have the transvectant identity
\begin{equation} 
(\alpha_\ux\beta_\ux,\gamma_\ux\delta_\ux)_1=\frac{1}{2}(\alpha\gamma)\beta_\ux\delta_\ux+\frac{1}{2}(\beta\delta)\alpha_\ux\gamma_\ux\ .
\label{identity.abcd} \end{equation}
\end{Lemma}

\proof
Write
\[
(\alpha_\ux\beta_\ux,\gamma_\ux\delta_\ux)_1= \ 
\parbox{3.5cm}{
\psfrag{a}{$\scriptstyle{\alpha}$}
\psfrag{b}{$\scriptstyle{\beta}$}
\psfrag{c}{$\scriptstyle{\gamma}$}
\psfrag{d}{$\scriptstyle{\delta}$}
\psfrag{x}{$\scriptstyle{x}$}
\psfrag{C}{$\scriptstyle{CG}$}
\includegraphics[width=3.5cm]{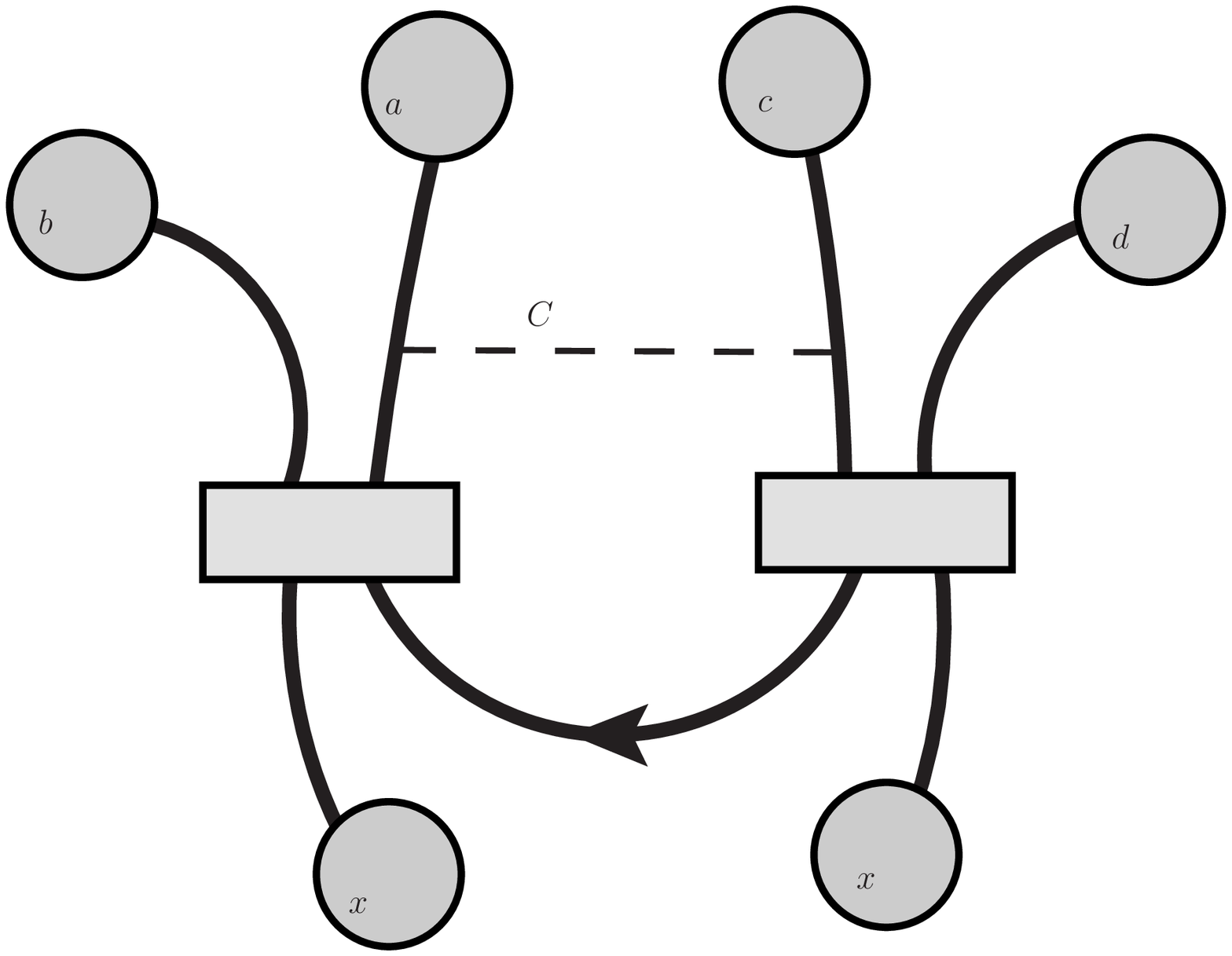}}
\]
and apply the Clebsch-Gordan (CG) identity
in~\cite[Eq. 2.9]{Abdesselam} at the place indicated by the dashed line.
This gives
\[
(\alpha_\ux\beta_\ux,\gamma_\ux\delta_\ux)_1=\ 
\parbox{3.5cm}{
\psfrag{a}{$\scriptstyle{\alpha}$}
\psfrag{b}{$\scriptstyle{\beta}$}
\psfrag{c}{$\scriptstyle{\gamma}$}
\psfrag{d}{$\scriptstyle{\delta}$}
\psfrag{x}{$\scriptstyle{x}$}
\includegraphics[width=3.5cm]{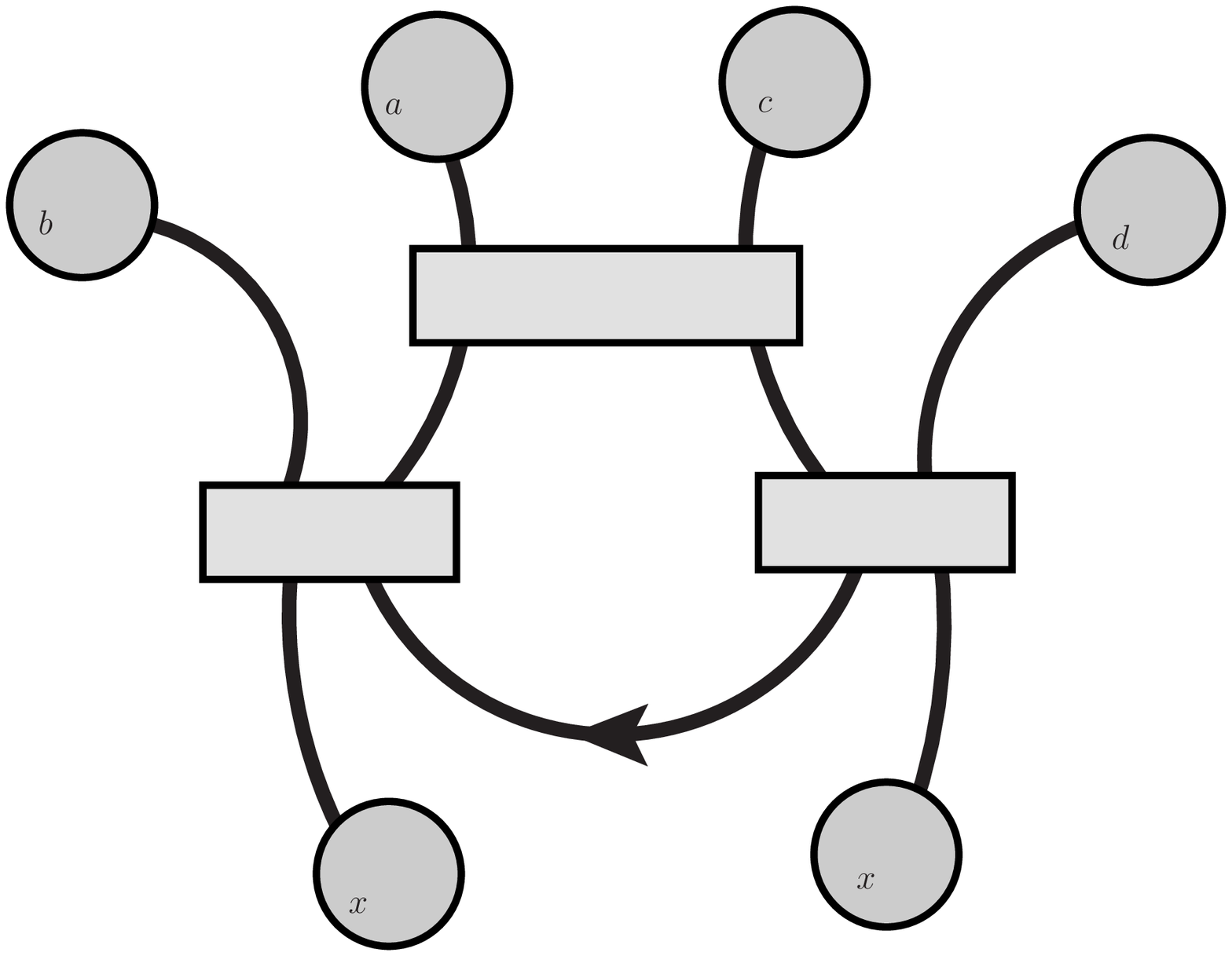}}
\ +\frac{1}{2}\ 
\parbox{3.5cm}{
\psfrag{a}{$\scriptstyle{\alpha}$}
\psfrag{b}{$\scriptstyle{\beta}$}
\psfrag{c}{$\scriptstyle{\gamma}$}
\psfrag{d}{$\scriptstyle{\delta}$}
\psfrag{x}{$\scriptstyle{x}$}
\includegraphics[width=3.5cm]{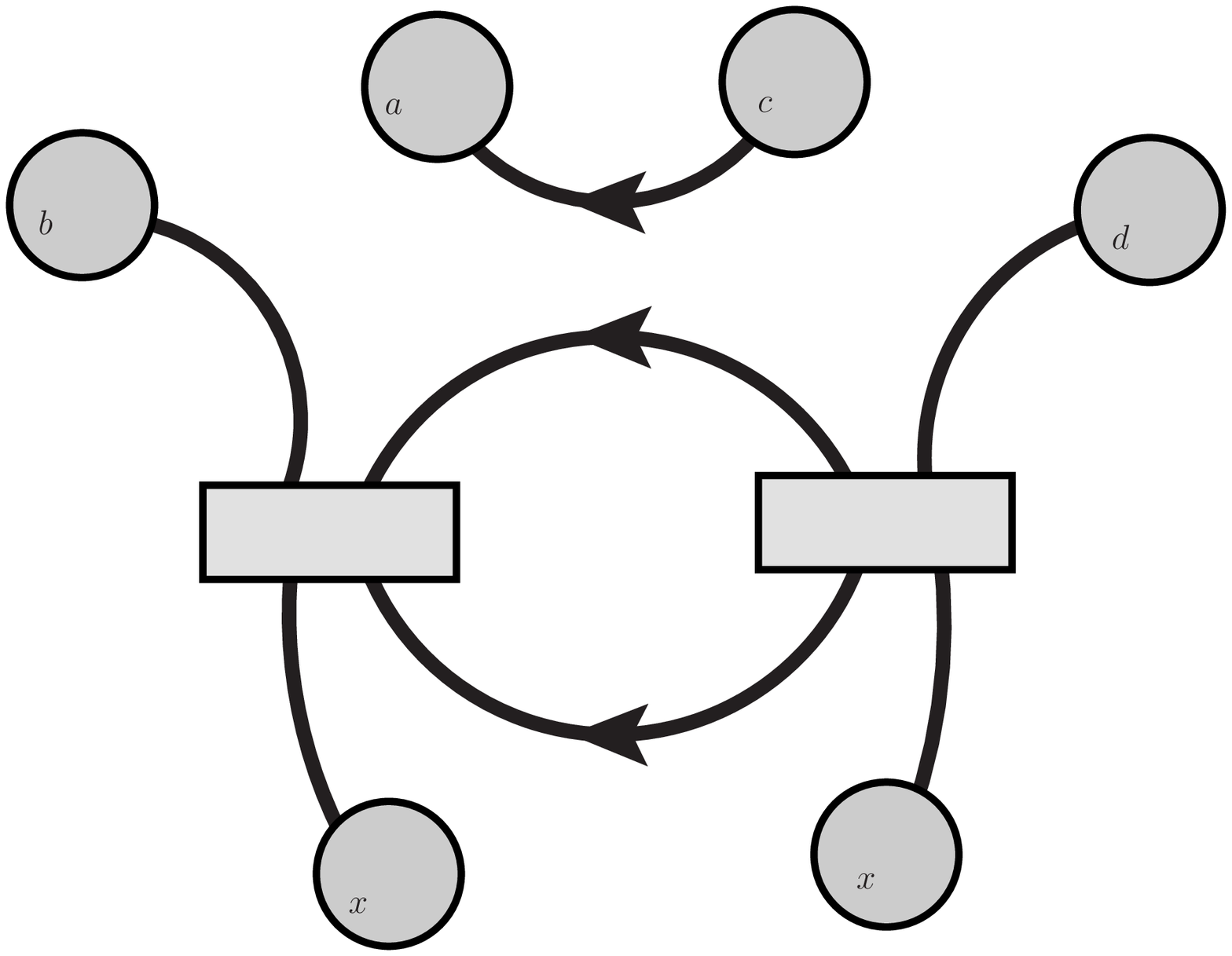}}\ \ .
\]
The weights $1$ and $\frac{1}{2}$ come from the ratios of binomial
coefficients in~\cite[Eq. 2.9]{Abdesselam}.  After expanding the symmetrizers, we get
\[
\parbox{2.2cm}{
\psfrag{b}{$\scriptstyle{\beta}$}
\psfrag{d}{$\scriptstyle{\delta}$}
\psfrag{x}{$\scriptstyle{x}$}
\includegraphics[width=2.2cm]{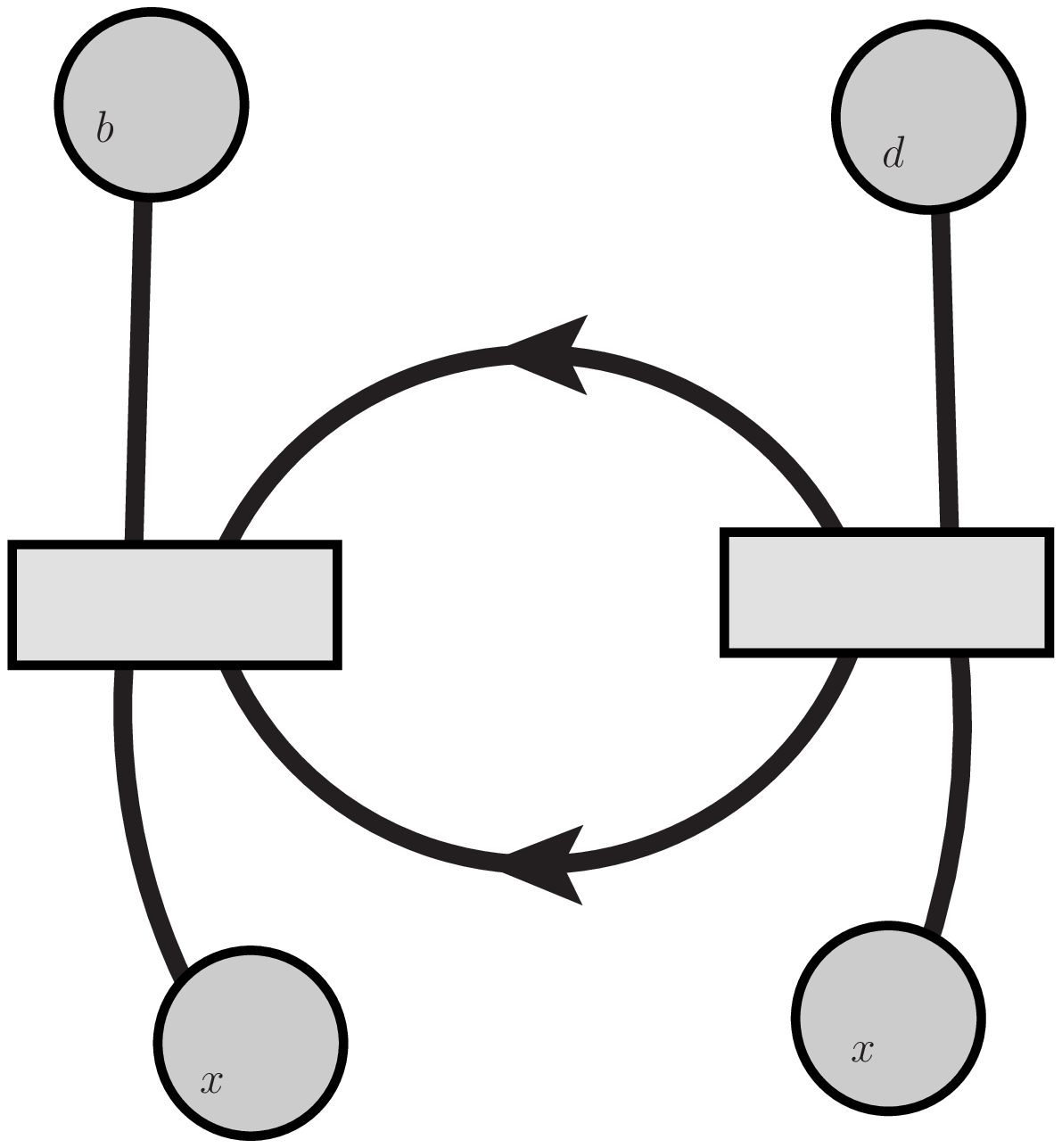}}
\ =\frac{1}{4}\left[\ 
\parbox{2.2cm}{
\psfrag{b}{$\scriptstyle{\beta}$}
\psfrag{d}{$\scriptstyle{\delta}$}
\psfrag{x}{$\scriptstyle{x}$}
\includegraphics[width=2.2cm]{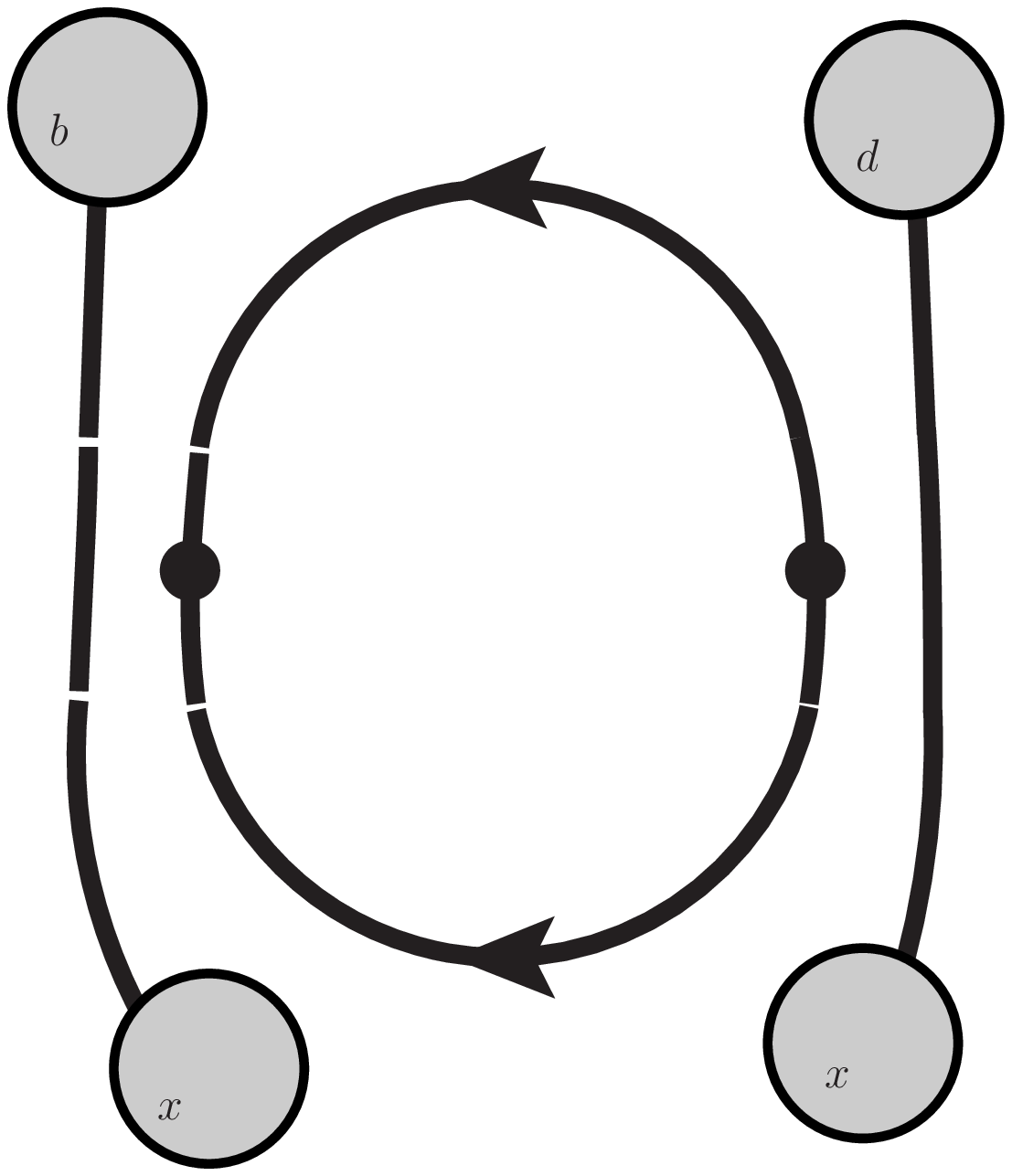}}
\ +\ 
\parbox{2.2cm}{
\psfrag{b}{$\scriptstyle{\beta}$}
\psfrag{d}{$\scriptstyle{\delta}$}
\psfrag{x}{$\scriptstyle{x}$}
\includegraphics[width=2.2cm]{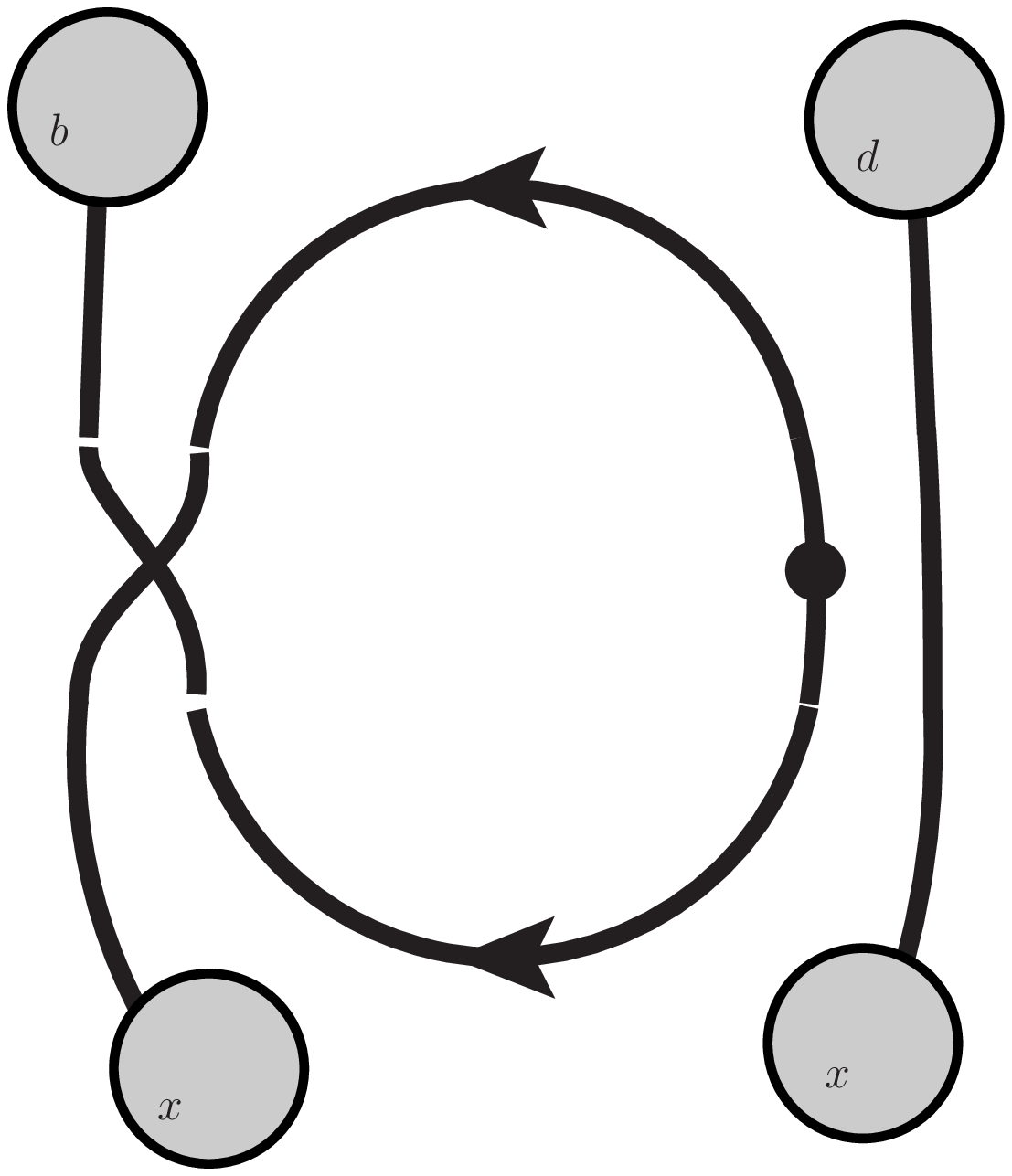}}
\ +\ 
\parbox{2.2cm}{
\psfrag{b}{$\scriptstyle{\beta}$}
\psfrag{d}{$\scriptstyle{\delta}$}
\psfrag{x}{$\scriptstyle{x}$}
\includegraphics[width=2.2cm]{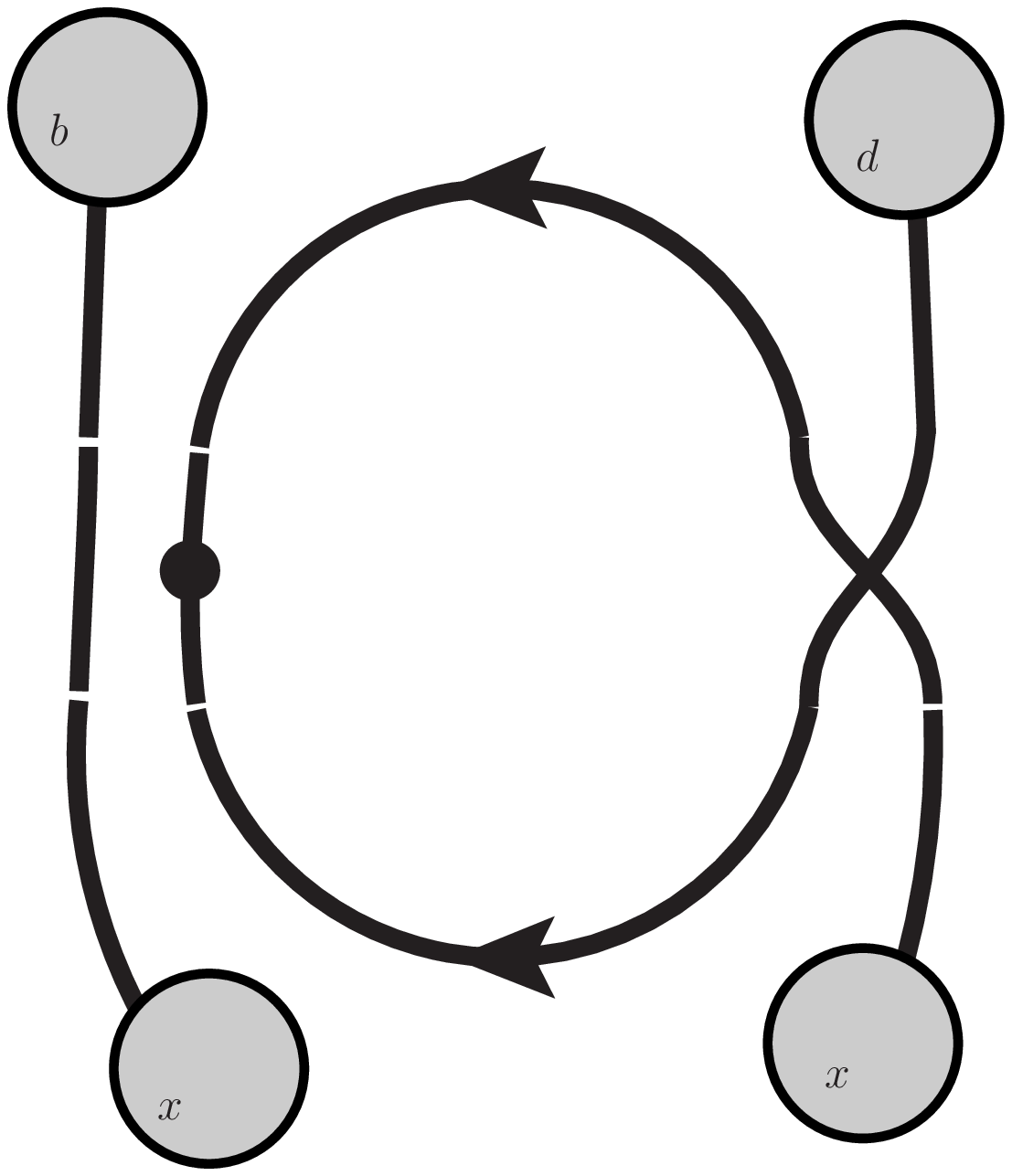}}
\ \right]=\beta_\ux\delta_\ux\ .
\]
Note that we haven't written the the fourth diagram with two
crossings, since it contains the bracket factor $(\ux\ux)=0$.
By applying the CG identity to the $\beta$ and $\delta$ strands, we get
\[
\parbox{3.5cm}{
\psfrag{a}{$\scriptstyle{\alpha}$}
\psfrag{b}{$\scriptstyle{\beta}$}
\psfrag{c}{$\scriptstyle{\gamma}$}
\psfrag{d}{$\scriptstyle{\delta}$}
\psfrag{x}{$\scriptstyle{x}$}
\includegraphics[width=3.5cm]{pic7.eps}}
\ =\ 
\parbox{3.5cm}{
\psfrag{a}{$\scriptstyle{\alpha}$}
\psfrag{b}{$\scriptstyle{\beta}$}
\psfrag{c}{$\scriptstyle{\gamma}$}
\psfrag{d}{$\scriptstyle{\delta}$}
\psfrag{x}{$\scriptstyle{x}$}
\includegraphics[width=3.5cm]{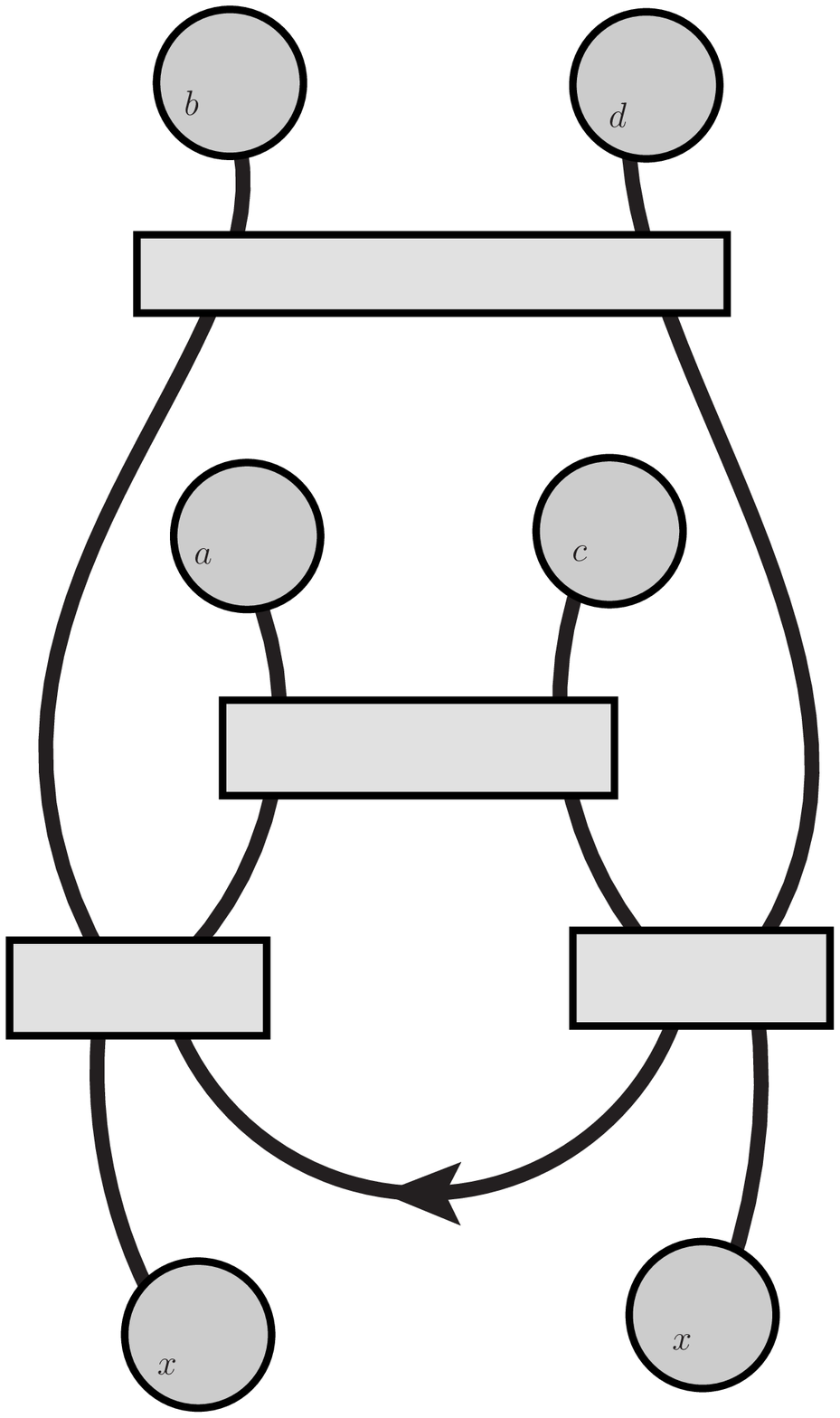}}
\ +\frac{1}{2}\ 
\parbox{3.5cm}{
\psfrag{a}{$\scriptstyle{\alpha}$}
\psfrag{b}{$\scriptstyle{\beta}$}
\psfrag{c}{$\scriptstyle{\gamma}$}
\psfrag{d}{$\scriptstyle{\delta}$}
\psfrag{x}{$\scriptstyle{x}$}
\includegraphics[width=3.5cm]{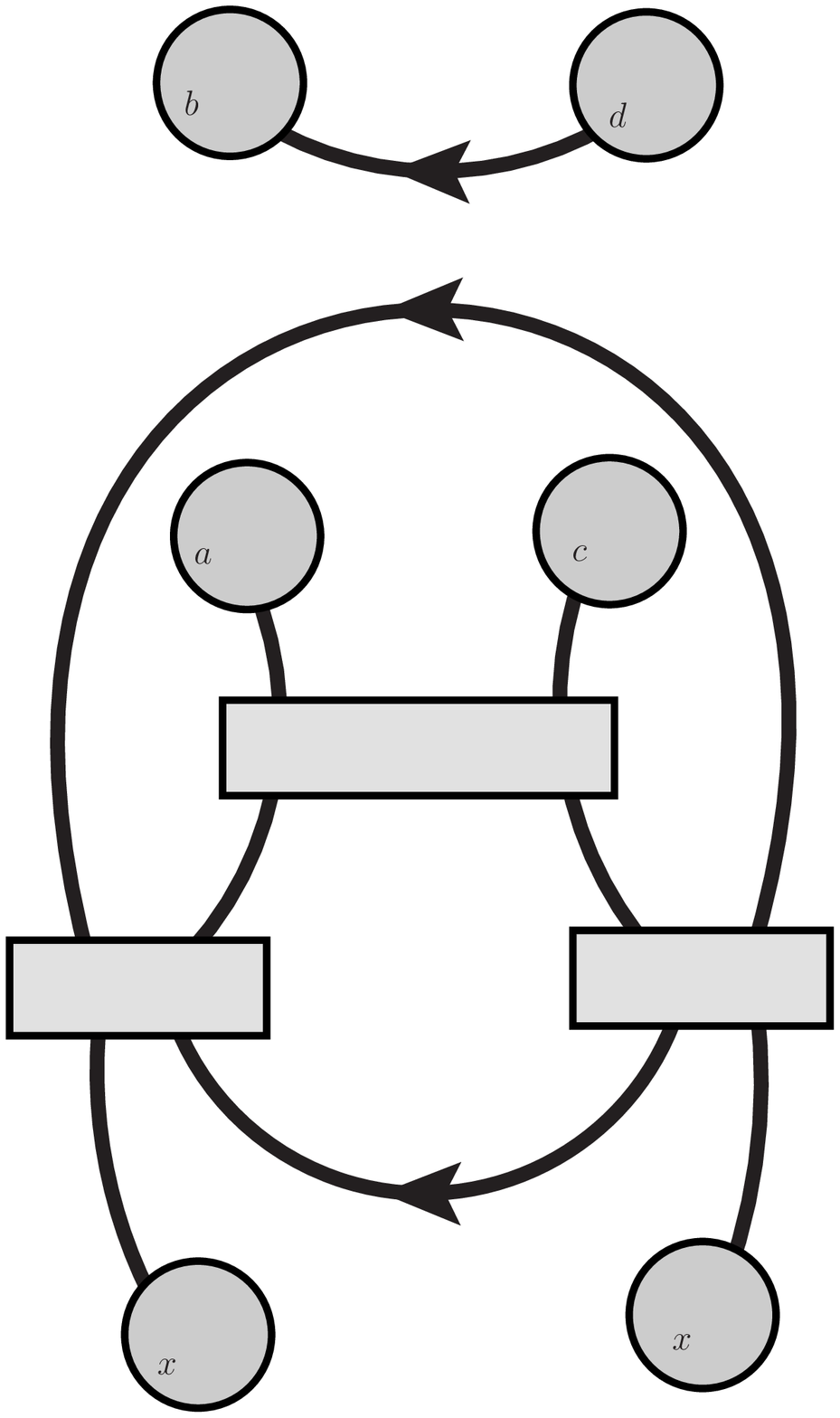}}\ \ .
\]
The second diagram can be computed by expanding the bottom two
symmetrizers as above, which gives the expression
$(\beta\delta)\alpha_\ux\gamma_\ux$.
We claim that the first diagram vanishes. Indeed, due to the presence of the top two symmetrizers,
if we move the bottom two symmetrizers so that they exchange places,
then the diagram becomes its own negative since this move
reverses the orientation of the bottom arrow. Now we get the required
identity by substituting back in the last equation for
$(\alpha_\ux\beta_\ux,\gamma_\ux\delta_\ux)_1$. \qed

\begin{Remark} \rm 
The left-hand side of~(\ref{identity.abcd}) corresponds to a pair partition $\{\{\alpha,\beta\},\{\gamma,\delta\}\}$.
We implicitly chose the `transverse' partition
$\{\{\alpha,\gamma\},\{\beta,\delta\}\}$ for the right-hand side. However,  we could have instead chosen
$\{\{\alpha,\delta\},\{\beta,\gamma\}\}$, which would give the equally valid identity
\[
(\alpha_\ux\beta_\ux,\gamma_\ux\delta_\ux)_1=
\frac{1}{2}(\alpha\delta)\beta_\ux\gamma_\ux+\frac{1}{2}(\beta\gamma)\alpha_\ux\delta_\ux. \]
If we average the last equality with~(\ref{identity.abcd}), the net result
is the `naive' four-term expansion of the transvectant as in~\cite[\S44 and \S49 (vii)]{GraceYoung}.
If one were to use the latter for a brute-force bracket monomial
computation of $\psi$, this would generate $4^3\times 2\times 3=384$
terms. (The factors of $4$ come from the calculation of $U$, $V$ and
$W$. The factor of $2$ comes from the computation of second
transvectants, and finally there are $3$ terms such as $(U,V)_2 W$.) 
Hence the previous lemma is essential in organizing the calculation of $\psi$ and reducing its complexity.
\end{Remark}

\subsection{} By Lemma~\ref{fourtotwo}, 
\[ U=\frac{1}{2}(cd)b_\ux f_\ux+\frac{1}{2} (bf) c_\ux d_\ux, \quad 
V=\frac{1}{2}(cd)a_\ux e_\ux+\frac{1}{2} (ae) c_\ux d_\ux. \]
Using the bilinearity of the second transvectant, we have
\[ \begin{aligned} 
4(U,V)_2 = & (cd)^2(b_\ux f_\ux,a_\ux e_\ux)_2+(cd)(ae)(b_\ux
f_\ux,c_\ux d_\ux)_2 + \\ 
& (bf)(cd)(c_\ux d_\ux,a_\ux e_\ux)_2+(bf)(ae)(c_\ux d_\ux,c_\ux
d_\ux)_2. \end{aligned} \] 
Now 
\[ (c_\ux d_\ux,c_\ux d_\ux)_2= \ 
\parbox{2.5cm}{
\psfrag{c}{$\scriptstyle{c}$}
\psfrag{d}{$\scriptstyle{d}$}
\includegraphics[width=2.5cm]{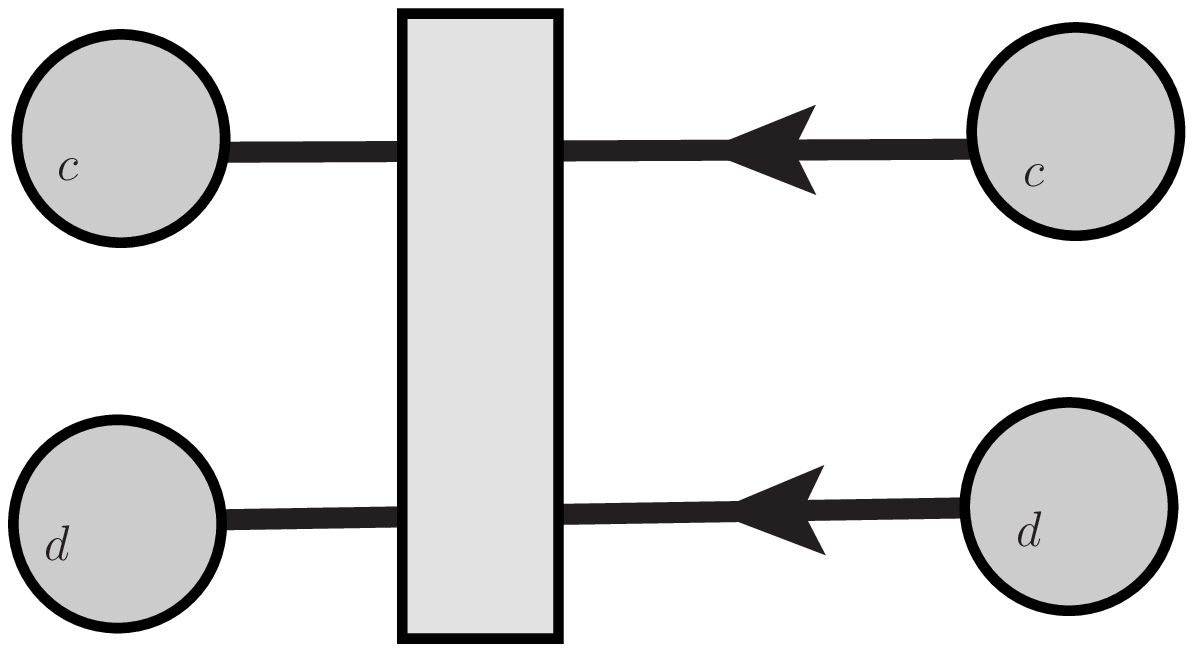}}
\ = -\frac{1}{2} (cd)^2
\]
and thus $(U,V)_2=\frac{1}{4}(cd)\ \exprS'$, where
\[
\exprS'=(cd)(b_x f_x,a_x e_x)_2+(ae)(b_x f_x,c_x d_x)_2
+(bf)(c_x d_x,a_x e_x)_2-\frac{1}{2}(ae)(bf)(cd)\ .
\]
We will show later that $\exprS'$ is in fact equal to the $\exprS$ of 
(\ref{defn.exprS}). Since second transvectants are symmetric bilinear forms, the previously
mentioned symmetries of $\exprS$ are particularly evident in the last equation.

By Lemma~\ref{fourtotwo}, we have 
\[
W=\frac{1}{2}(ae)b_\ux f_\ux+\frac{1}{2}(bf)a_\ux e_\ux\ .
\]
This results in
\begin{equation}
(U,V)_2 W=\frac{1}{8}\exprS'\times\{
(cd)(ae)b_\ux f_\ux+(cd)(bf)a_\ux e_\ux
\}\ . \label{UVWSeq} \end{equation}
The exchange of letters $b \leftrightarrow c, d \leftrightarrow f$
brings about an exchange of $V$ and $W$. Applying this to (\ref{UVWSeq}) gives
\begin{equation}
(U,W)_2 V=\frac{1}{8} \, \exprS' \times\{
(bf)(ae)c_\ux d_\ux+(bf)(cd)a_\ux e_\ux
\}\ . \label{UWVSeq} \end{equation}

Likewise, the exchange $a \leftrightarrow c, d \leftrightarrow e$ 
exchanges $U$ and $W$. Applying this to (\ref{UVWSeq}) gives
\begin{equation}
(W,V)_2 U=\frac{1}{8} \, \exprS' \times\{
(ae)(cd)b_\ux f_\ux+(ae)(bf)c_\ux d_\ux
\}\ .
\label{WVUSeq}
\end{equation}

Now substitute (\ref{UVWSeq}), (\ref{UWVSeq}), and (\ref{WVUSeq})
in the result of Lemma~\ref{psilem} and simplify. This gives the
required formula for $\psi$.

\subsection{} We now proceed with the simplification of $\exprS'$.
By expanding the symmetrizers implicit in the three second
transvectants, we get
\[ \begin{aligned} 
2\exprS' = & (cd)(ba)(fe)+(cd)(be)(fa)+(ae)(bc)(fd) + \\
& (ae)(bd)(fc)+(bf)(ca)(de)+(bf)(ce)(da)-(ae)(bf)(cd). 
\end{aligned} \] 
Now insert the GP relation $(ba)(fe)=(bf)(ae)-(be)(af)$ in the first
term, and similarly the relations 
\[ 
(bc)(fd)=(bf)(cd)-(bd)(cf), 
(ca)(de)=(cd)(ae)-(ce)(ad), \] 
respectively in the third and the fifth term. After an
expansion, cancellation and a division by $2$, we get 
\[
\exprS' =(cd)(be)(fa)+(ae)(bd)(fc)+(bf)(ce)(da)+(ae)(bf)(cd). 
\]
Now insert the GP relations 
\[ (cd)(be)=(cb)(de)-(ce)(db), \quad 
(bf)(ce)=(bc)(fe)-(be)(fc)
\]
respectively in the first and the third term, to get 
\[
\exprS' = \underbrace{-(ce)(db)(fa)-(be)(fc)(da)}_{\exprS}+\exprT, 
\]
where 
\[ \exprT=(cb)(de)(fa)+(ae)(bd)(fc)+(bc)(fe)(da)+(ae)(bf)(cd).  \]
We only need to verify that $\exprT$ is identically zero, which would
imply $\exprS' = \exprS$. To this end, insert the GP relations 
\[
(de)(fa)=(df)(ea)-(da)(ef), \quad (bd)(fc)=(bf)(dc)-(bc)(df)
\]
respectively in the first and second term of $\exprT$. 
The six resulting terms cancel in pairs, and thus $\exprT=0$. This
completes the proof of Proposition~\ref{proposition.ae}. \qed 

\subsection{} 
The invariant $\exprS$ has played an important role in the proof.
The following proposition gives another notable property of this invariant.

\begin{Proposition} \rm 
The polynomial $\exprS$ and the simpler expression $(ae)(bf)(cd)$ 
form a basis of the vector space of multilinear $SL_2$-invariants of $a,b,\ldots,f$
which satisfy the $\SG_3\times \ZZ_2$ symmetry mentioned in Remark~\ref{sym.remark}.
\end{Proposition}

\proof
We first show that the two invariants are not proportional. Indeed,
$\exprS$ is not expressible as a bracket monomial and thus its 
expression in~(\ref{defn.exprS}) is as simple as possible.
This can be seen  by making the usual specialization of sending three points to 0, 1, and $\infty$, i.e.,
letting say $a=(0,1)$, $b=(1,1)$, $c=(1,0)$, $d=(x,1)$, $e=(y,1)$ and $f=(z,1)$.
One then gets $\exprS = -xy+x+z-xz$. 

A bracket monomial would have a bracket containing $c$ which gives $\pm 1$. The remaining two brackets
would give affine linear expressions in $x,y,z$. If $\exprS$ were proportional to a bracket monomial, then
the polynomial $-xy+x+z-xz$ would be reducible. If one homogenizes by
adding a variable $t$, then 
\[
-xy+xt+zt-xz=\frac{1}{2} XMX^{\rm T}
\]
where $X=(x,y,z,t)$ and $M=\left(
\begin{array}{cccc}
0 & -1 & -1 & 1\\
-1 & 0 & 0 & 0\\
-1 & 0 & 0 & 1\\
1 & 0 & 1 & 0
\end{array}
\right)$. 
Since $\det(M)=1\neq 0$, the polynomial above is irreducible, which proves our claim.

We now show that the vector space under consideration has dimension two.
Introduce the invariants
\[ \begin{array}{lll} 
B_1 = (ae)(bf)(cd), & B_2 = (ab)(ec)(fd), & B_3 =  (ad)(bc)(ef), \\ 
B_4 = (ab)(cd)(fe), & B_5 = (ea)(bc)(df). 
\end{array} \] 
It is a consequence of Kempe's 
Circular Straightening Theorem (see, e.g,~\cite[Prop. 2.6]{HowardMSV} or~\cite[Lemma 6.2]{KungR})
that $B_1,\ldots,B_5$ form a basis of the space of multilinear $SL_2$-invariants of the six points $a,b,\ldots,f$.
Indeed, if we order these points cyclically as $a,b,c,d,f,e$, then $B_1,\ldots,B_5$
correspond to the five non-crossing chord configurations.

Let $J, K, L$ be as in Remark~\ref{sym.remark}. A straightforward calculation shows that the action of these
generators in the $B$-basis is given by the following matrices: 
\[
J=\left(
\begin{array}{rrrrr}
1 & 0 & 0 & 0 & 0 \\
0 & -1 & 0 & 0 & -1 \\
0 & 0 & -1 & 0 & 1 \\
0 & -1 & 1 & 1 & -1 \\
0 & 0 & 0 & 0 & 1
\end{array} \right)\ ,\ 
K=\left(
\begin{array}{rrrrr}
1 & 0 & 0 & 0 & 0 \\
0 & -1 & 0 & 1 & 0 \\
0 & 0  & -1 & -1 & 0 \\
0 & 0 & 0 & 1 & 0 \\
0 & 1 & -1 & -1 & 1
\end{array}
\right)\ ,
L=\left(
\begin{array}{rrrrr}
-1 & 0 & 0 & 0 & 0 \\
0 & 0 & 1 & 0 & 0 \\
0 & 1 & 0 & 0 & 0 \\
0 & 0 & 0 & -1 & 0\\
0 & 0 & 0  & 0 & -1
\end{array}
\right)\ .
\]
The permutations typically create crossings (at most two), and the latter can be undone using a GP relation
to express the result in the $B$-basis.  This procedure gives the matrices above. 
There are~\emph{a priori} fifteen equations defining the intersection of 
${\rm Ker}(J-I)$, ${\rm Ker}(K-I)$ and ${\rm Ker}(L+I)$, but they 
reduce to a homogeneous system of three independent equations given by the matrix
\[
\left(
\begin{array}{rrrrr}
0 & -2 & 0 & 0 & -1 \\
0 & 0 & -2 & 0 & 1 \\
0 & 0 & 0 & 1 & 1
\end{array}
\right)\ .
\]
Therefore the dimension of the solution space is two.
The invariant $\exprS$ corresponds to the coordinate vector
$(-2,-1,1,-2,2)^{\rm T}$, 
which of course satisfies this homogeneous system. \qed

\begin{Remark} \rm 
There is a simple combinatorial recipe for finding the two bracket
monomials appearing in $\exprS$.
Draw the oriented graph on six vertices given by the edges
$a\leftarrow e$, $b\leftarrow f$, $c\leftarrow d$, which correspond to the three
quadratics used to build $U$, $V$ and $W$. Now ask: how can one add three more directed edges in order to form a
\emph{properly oriented} 6-cycle? The two possible answers give the
two required bracket monomials. 
\end{Remark}

\subsection{Proof of Proposition~\ref{proposition.ax}}
Recall that $U, V, W$ are now arbitrary quadratics, and $a_\ux$ is a
linear form. 
By expanding the symmetrizer, we have 
\[
\left(U, (V,a_\ux)_1 (W,a_\ux)_1\right)_1=\ 
\parbox{3cm}{
\psfrag{U}{$\scriptstyle{U}$}
\psfrag{V}{$\scriptstyle{V}$}
\psfrag{W}{$\scriptstyle{W}$}
\psfrag{x}{$\scriptstyle{x}$}
\psfrag{a}{$\scriptstyle{a}$}
\includegraphics[width=3cm]{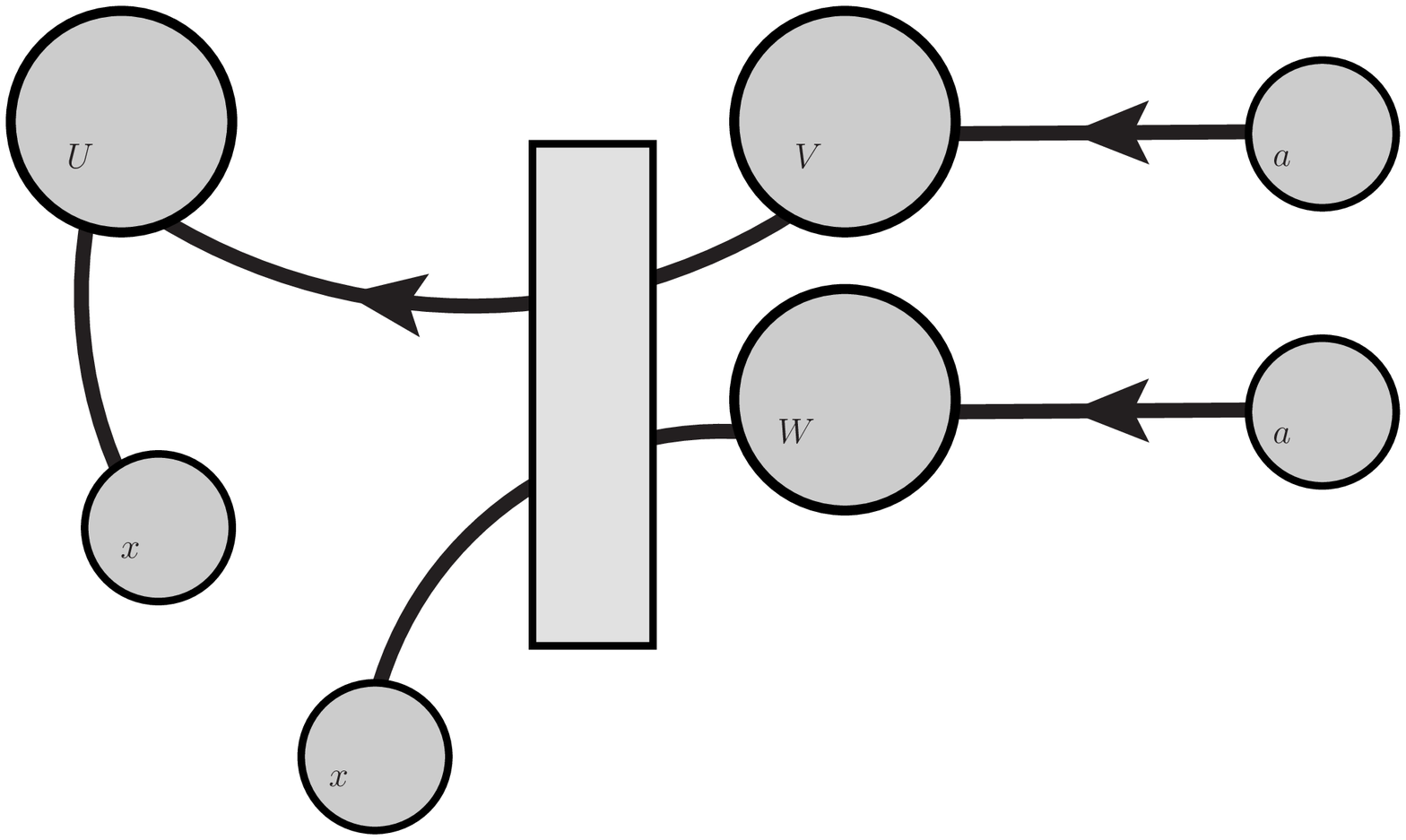}}
\ 
=\frac{1}{2}G_{\rm v}+\frac{1}{2}G_{\rm w}
\]
with
\[
G_{\rm v}=\ 
\parbox{3cm}{
\psfrag{U}{$\scriptstyle{U}$}
\psfrag{V}{$\scriptstyle{V}$}
\psfrag{W}{$\scriptstyle{W}$}
\psfrag{x}{$\scriptstyle{x}$}
\psfrag{a}{$\scriptstyle{a}$}
\includegraphics[width=3cm]{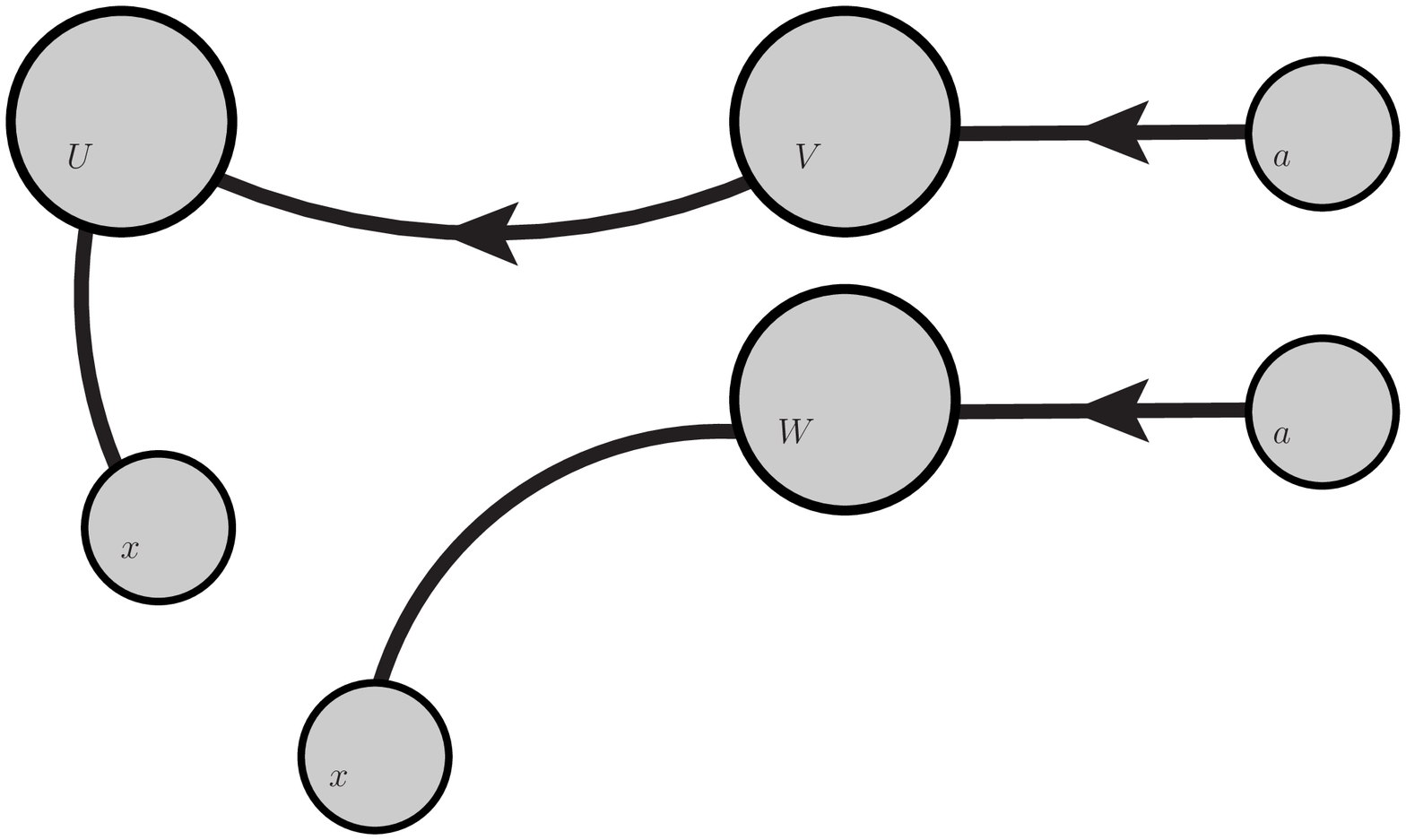}}
\]
and
\[
G_{\rm w}=\ 
\parbox{3cm}{
\psfrag{U}{$\scriptstyle{U}$}
\psfrag{V}{$\scriptstyle{V}$}
\psfrag{W}{$\scriptstyle{W}$}
\psfrag{x}{$\scriptstyle{x}$}
\psfrag{a}{$\scriptstyle{a}$}
\includegraphics[width=3cm]{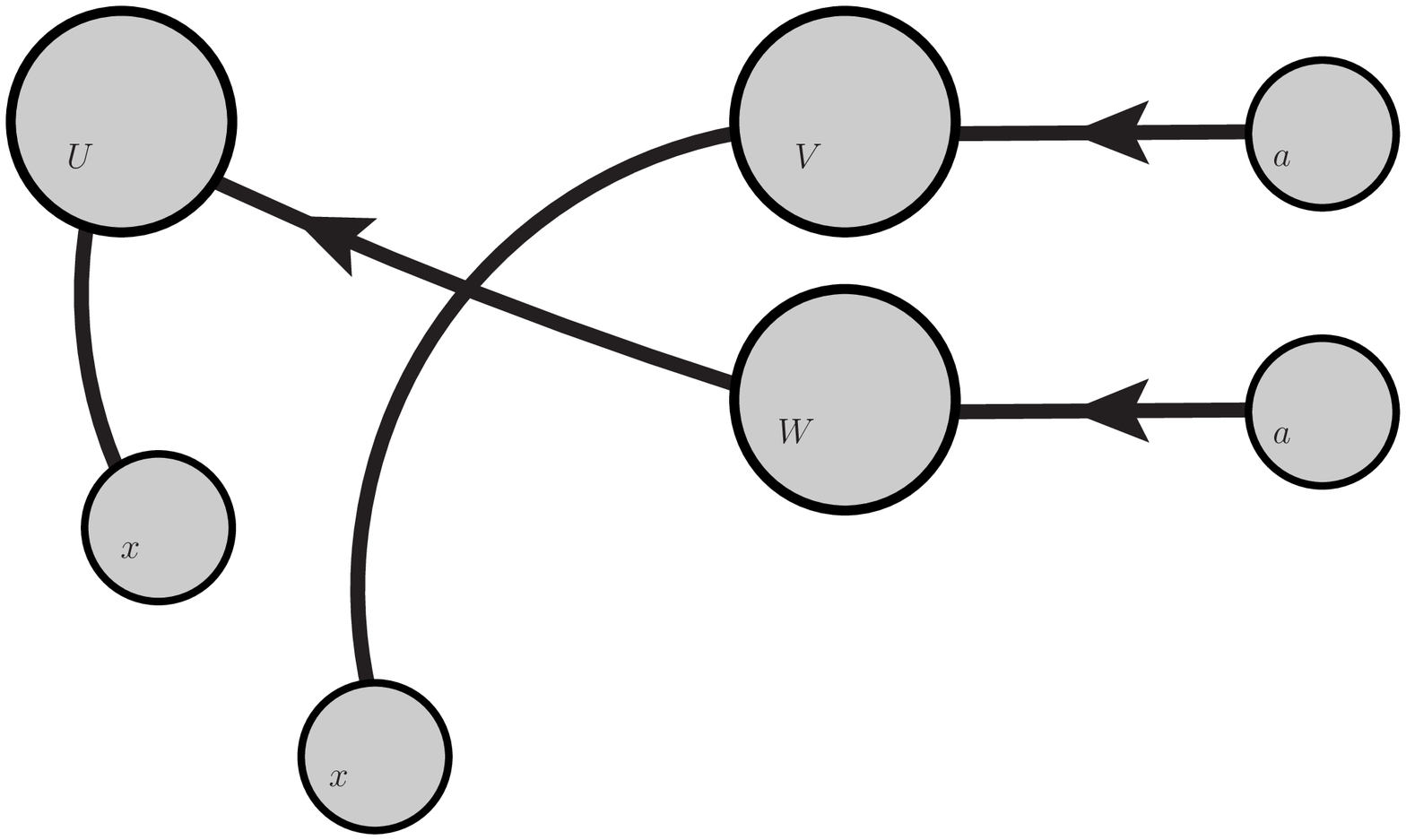}}\ \ .
\]
We will compute $G_{\rm v}$
and deduce the analogous formula for $G_{\rm w}$ by exchanging $V$ and $W$.
One can rewrite
\[
G_{\rm v}=\ 
\parbox{3cm}{
\psfrag{U}{$\scriptstyle{U}$}
\psfrag{V}{$\scriptstyle{V}$}
\psfrag{W}{$\scriptstyle{W}$}
\psfrag{x}{$\scriptstyle{x}$}
\psfrag{a}{$\scriptstyle{a}$}
\includegraphics[width=3cm]{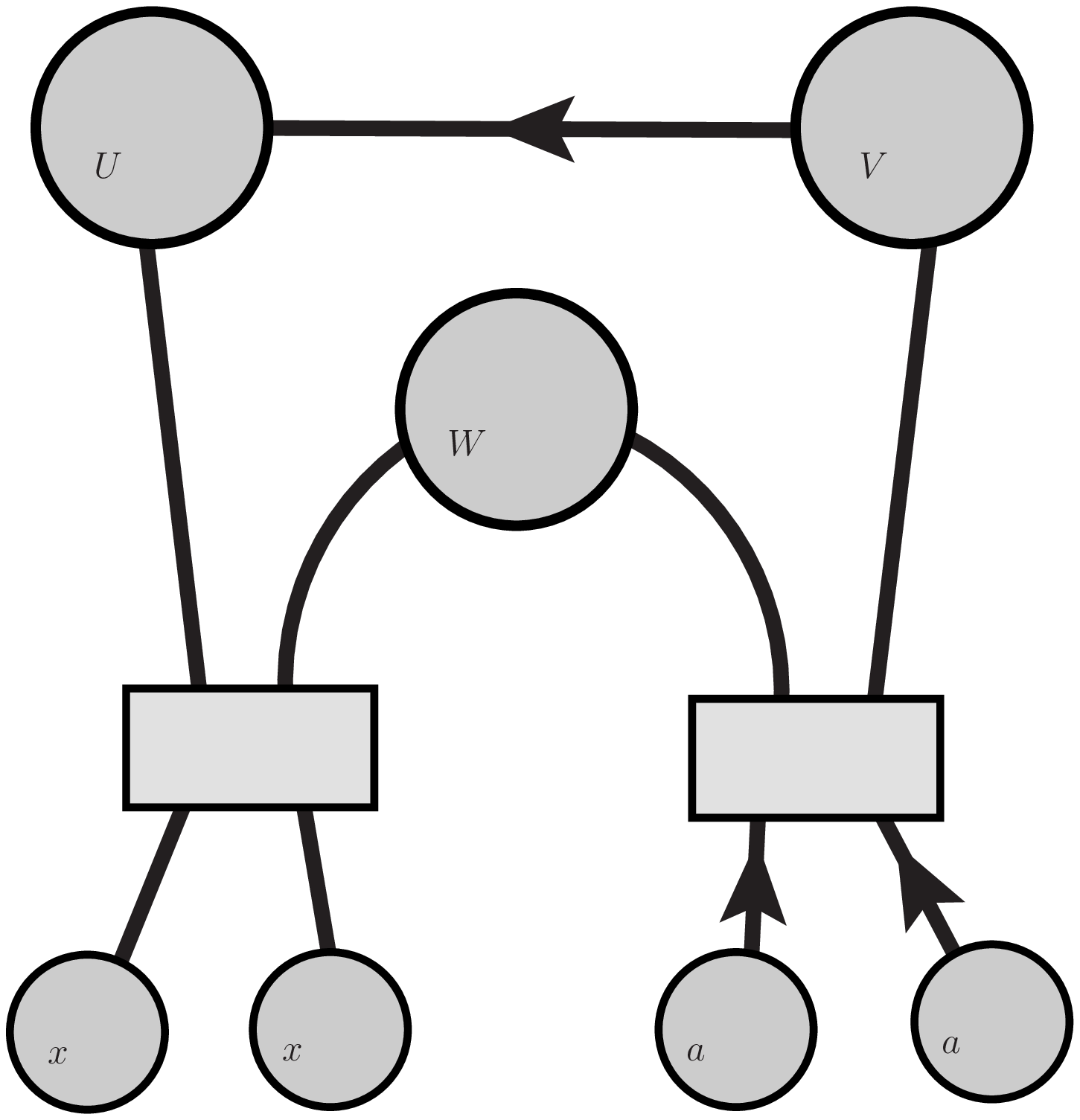}}
\]
and apply the CG identity~\cite[Eq. 2.9]{Abdesselam} between the bottom two symmetrizers. 
This results in $G_{\rm v}=G_{{\rm v}0}+G_{{\rm v}1}+\frac{1}{3}G_{{\rm v}2}$
with
\[
G_{{\rm v}0}=\ 
\parbox{3cm}{
\psfrag{U}{$\scriptstyle{U}$}
\psfrag{V}{$\scriptstyle{V}$}
\psfrag{W}{$\scriptstyle{W}$}
\psfrag{x}{$\scriptstyle{x}$}
\psfrag{a}{$\scriptstyle{a}$}
\includegraphics[width=3cm]{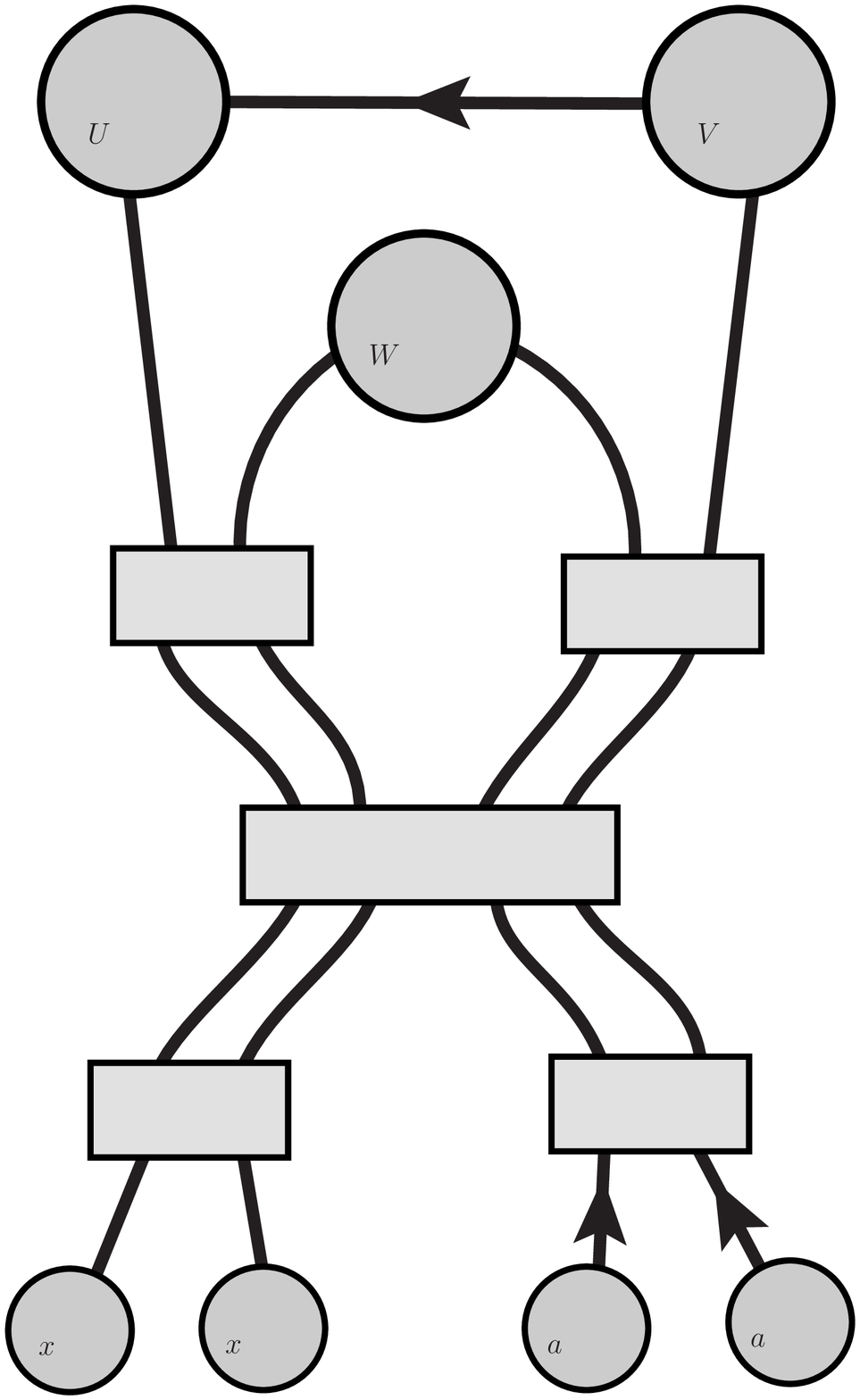}}
\ ,\ \ 
G_{{\rm v}1}=\ 
\parbox{3cm}{
\psfrag{U}{$\scriptstyle{U}$}
\psfrag{V}{$\scriptstyle{V}$}
\psfrag{W}{$\scriptstyle{W}$}
\psfrag{x}{$\scriptstyle{x}$}
\psfrag{a}{$\scriptstyle{a}$}
\includegraphics[width=3cm]{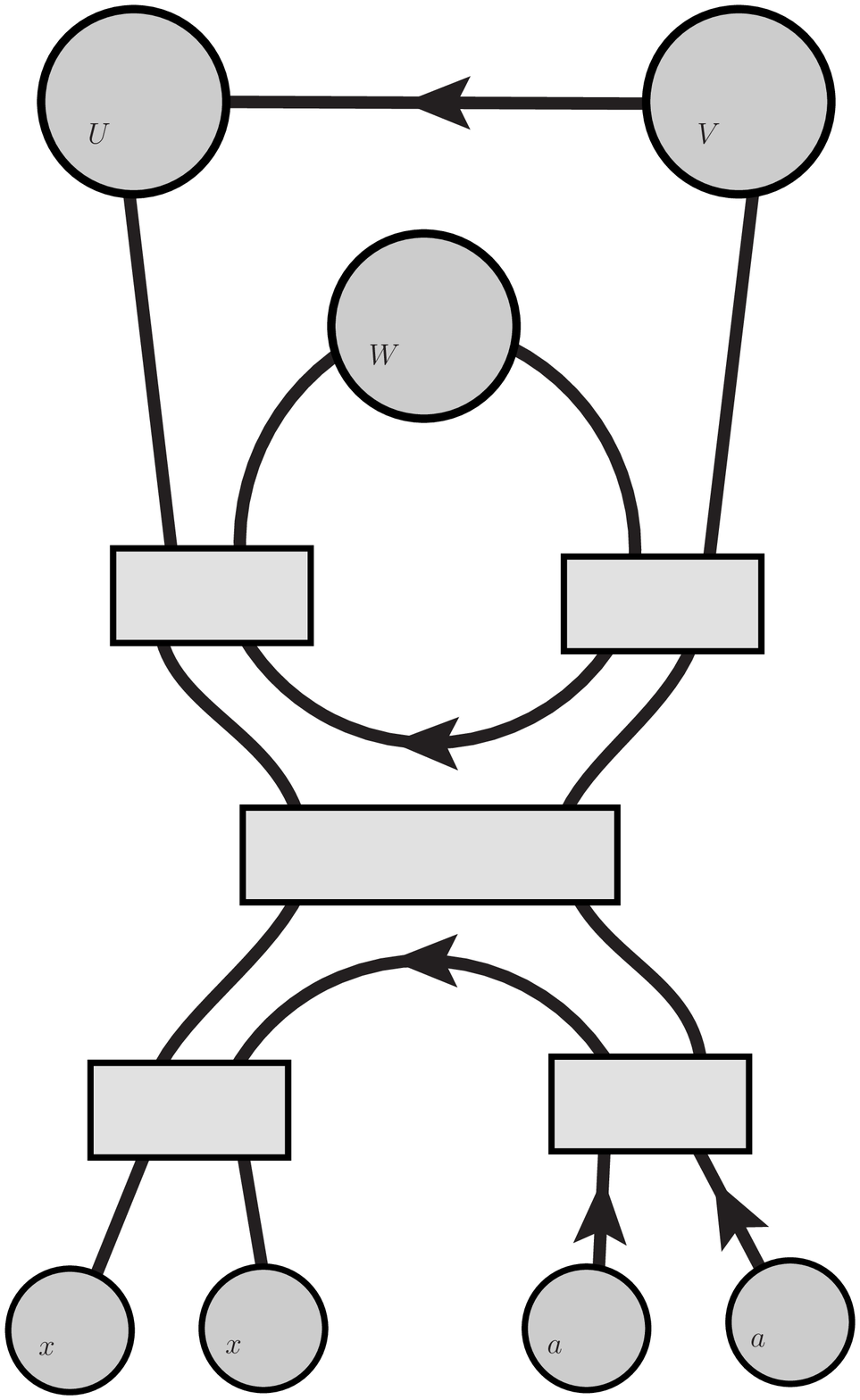}}
\ ,\ \ 
{\rm and}
\ \ 
G_{{\rm v}2}=\ 
\parbox{3cm}{
\psfrag{U}{$\scriptstyle{U}$}
\psfrag{V}{$\scriptstyle{V}$}
\psfrag{W}{$\scriptstyle{W}$}
\psfrag{x}{$\scriptstyle{x}$}
\psfrag{a}{$\scriptstyle{a}$}
\includegraphics[width=3cm]{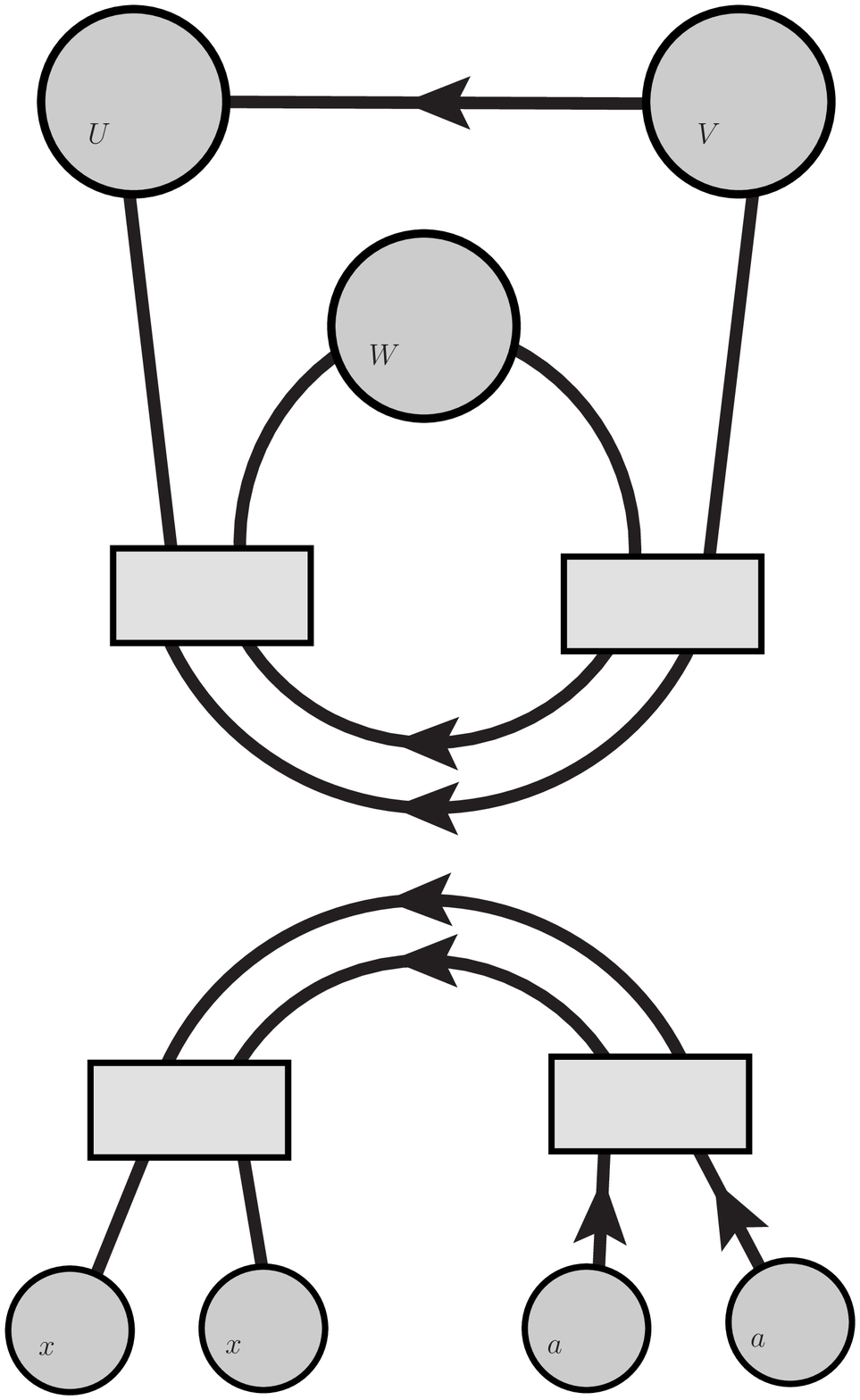}}
\ .
\]
Having the big symmetrizer eat up the smaller ones, we can write
\[
G_{{\rm v}0}=\ 
\parbox{3cm}{
\psfrag{U}{$\scriptstyle{U}$}
\psfrag{V}{$\scriptstyle{V}$}
\psfrag{W}{$\scriptstyle{W}$}
\psfrag{x}{$\scriptstyle{x}$}
\psfrag{a}{$\scriptstyle{a}$}
\includegraphics[width=3cm]{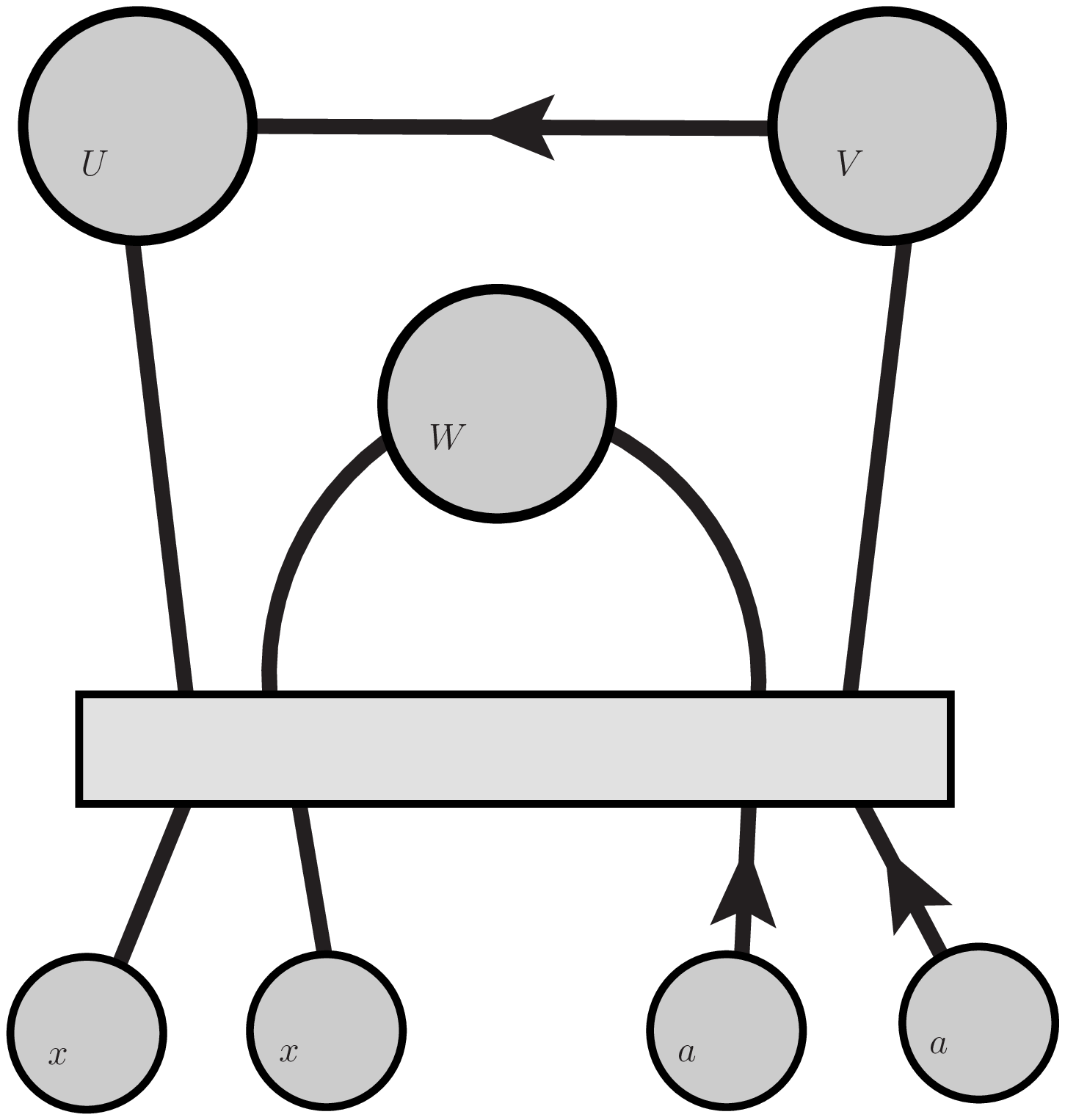}}
\ 
=\left((U,V)_1 W, a_{\ux}^2\right)_2\ .
\]
Passing the bottom arrows through the right symmetrizer and using
idempotence, we get
\[
\parbox{2.5cm}{
\psfrag{U}{$\scriptstyle{U}$}
\psfrag{V}{$\scriptstyle{V}$}
\psfrag{W}{$\scriptstyle{W}$}
\includegraphics[width=2.5cm]{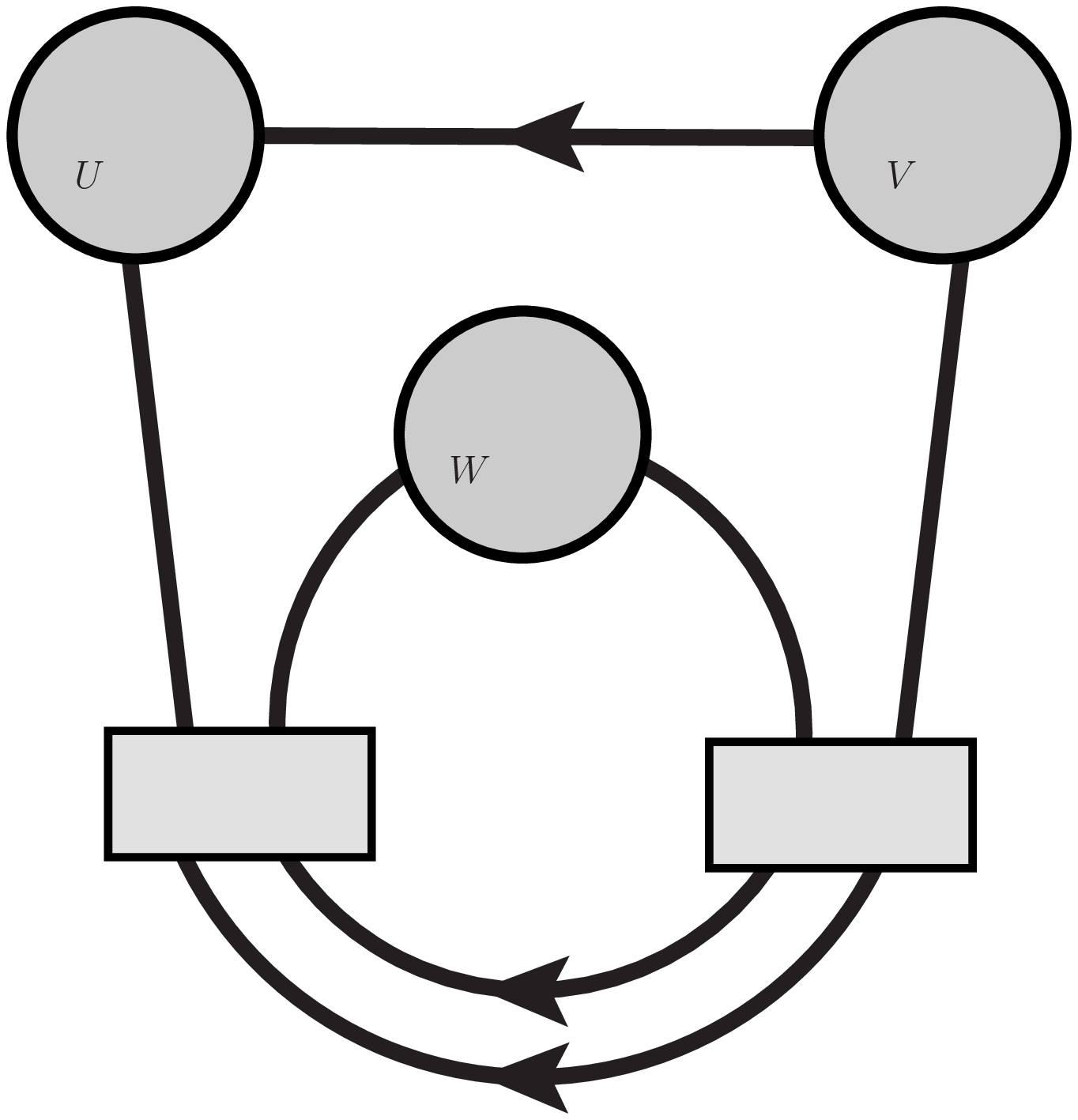}}
\ =\ 
\parbox{2.5cm}{
\psfrag{U}{$\scriptstyle{U}$}
\psfrag{V}{$\scriptstyle{V}$}
\psfrag{W}{$\scriptstyle{W}$}
\includegraphics[width=2.5cm]{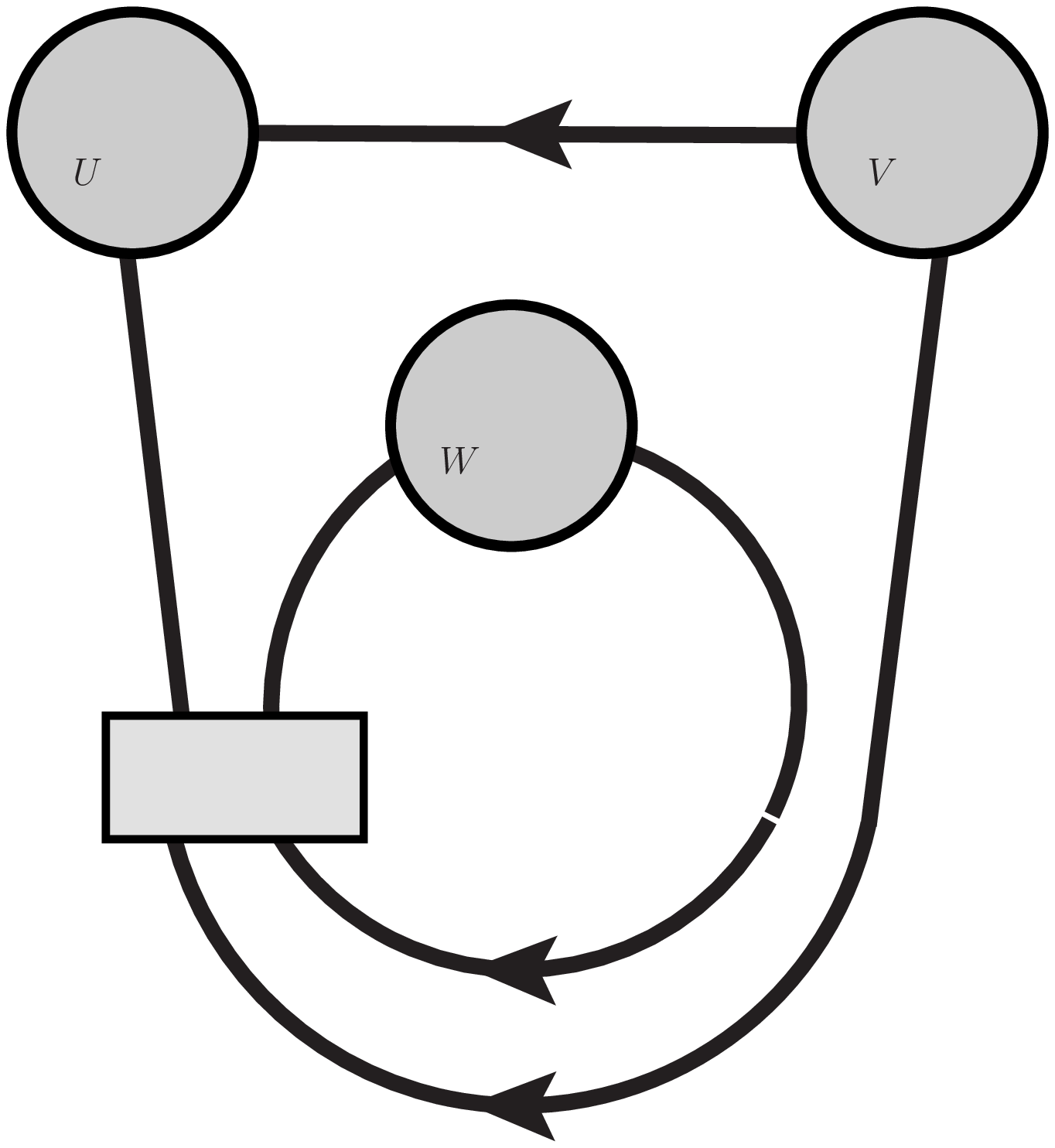}}
\ .
\]
Now expand the symmetrizer and ignore the vanishing term with the $W$ self-loop.
This gives
\[
\frac{1}{2}=\ 
\parbox{2.5cm}{
\psfrag{U}{$\scriptstyle{U}$}
\psfrag{V}{$\scriptstyle{V}$}
\psfrag{W}{$\scriptstyle{W}$}
\includegraphics[width=2.5cm]{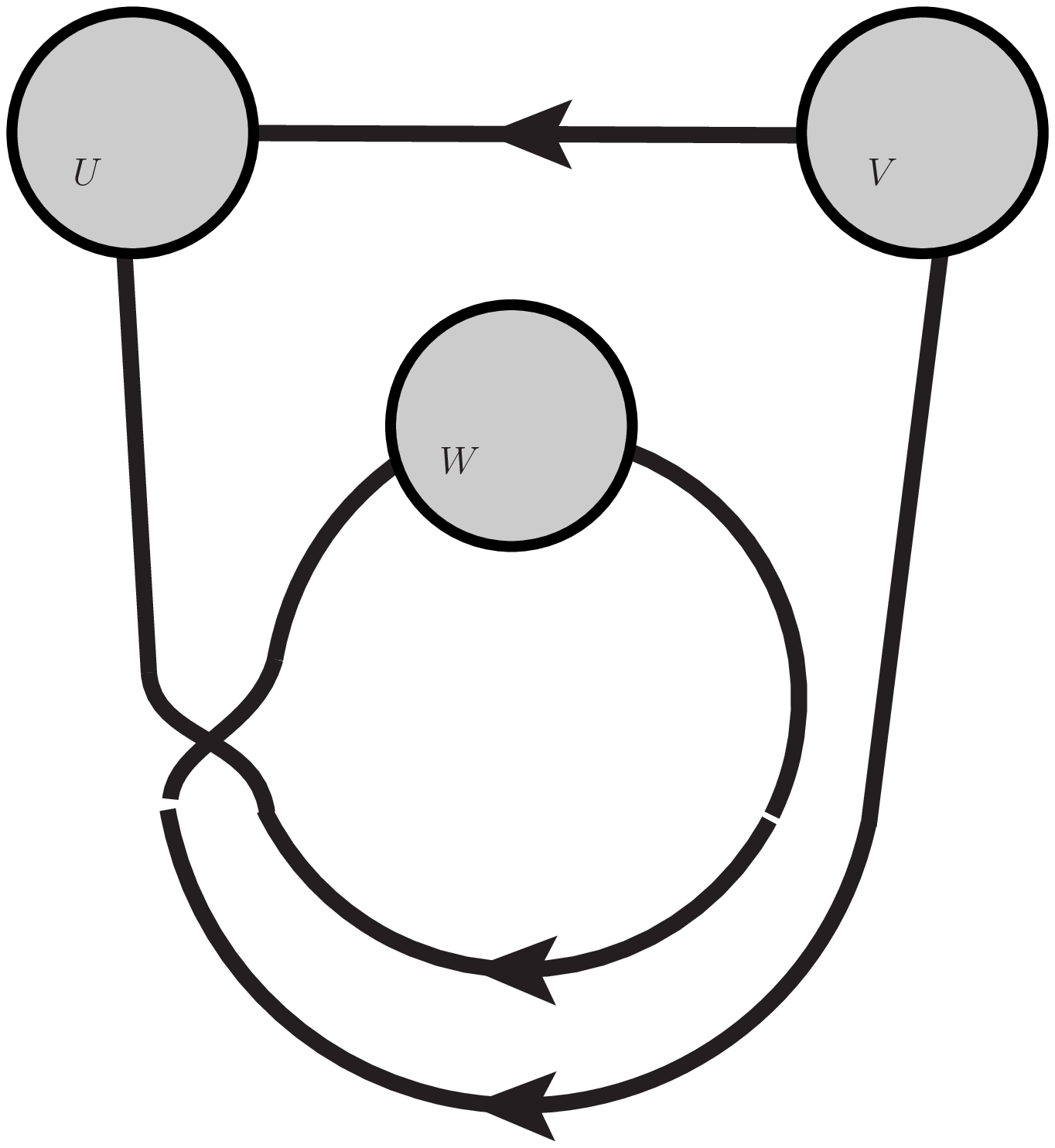}}
\ =\frac{1}{2}\ 
\parbox{2.5cm}{
\psfrag{U}{$\scriptstyle{U}$}
\psfrag{V}{$\scriptstyle{V}$}
\psfrag{W}{$\scriptstyle{W}$}
\includegraphics[width=2.5cm]{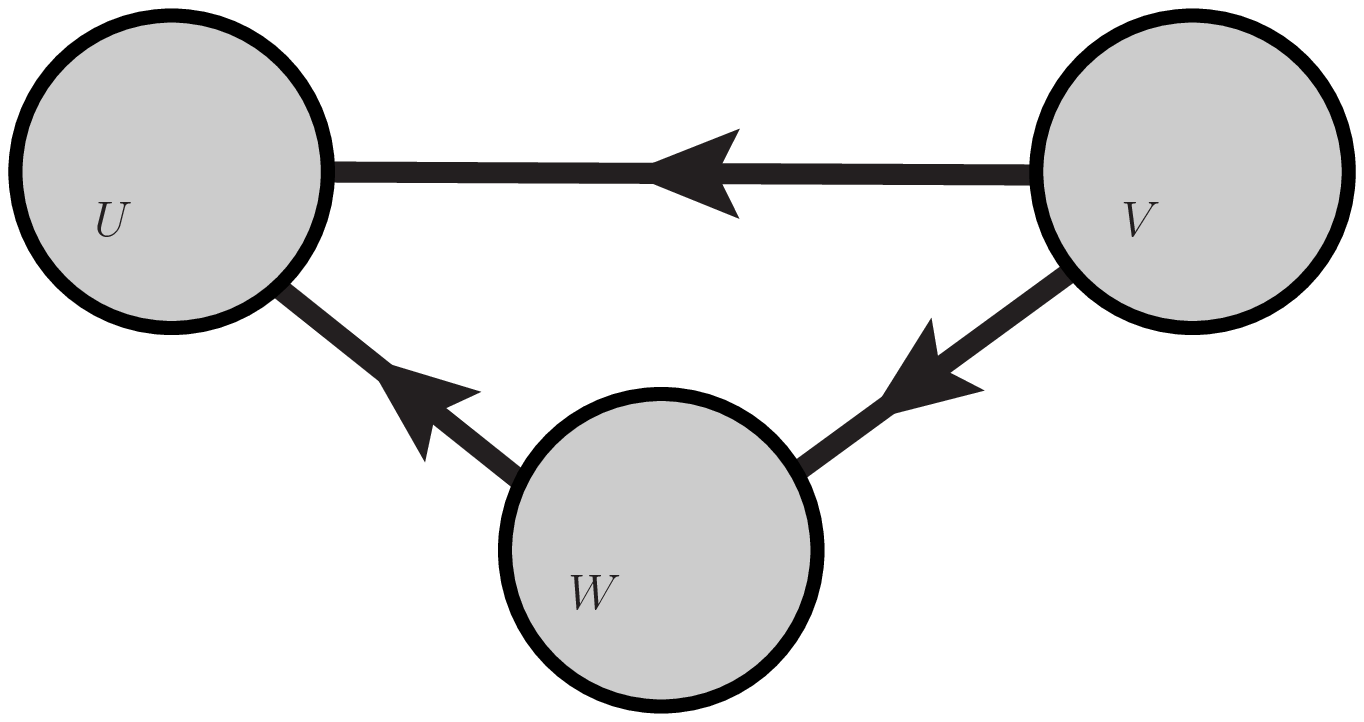}}
\ =0\ .
\]
Indeed, exchanging the positions of the $V$ and $W$ blobs shows that the diagram is equal to its negative.
Thus $G_{{\rm v}2}=0$.

Now remove the redundant $\ux$-symmetrizer and pass the arrows through
the symmetrizer on the bottom right of the previous diagram for
$G_{{\rm v}1}$. Then we have 
\[
G_{{\rm v}1}=\ 
\parbox{3cm}{
\psfrag{U}{$\scriptstyle{U}$}
\psfrag{V}{$\scriptstyle{V}$}
\psfrag{W}{$\scriptstyle{W}$}
\psfrag{x}{$\scriptstyle{x}$}
\psfrag{a}{$\scriptstyle{a}$}
\includegraphics[width=3cm]{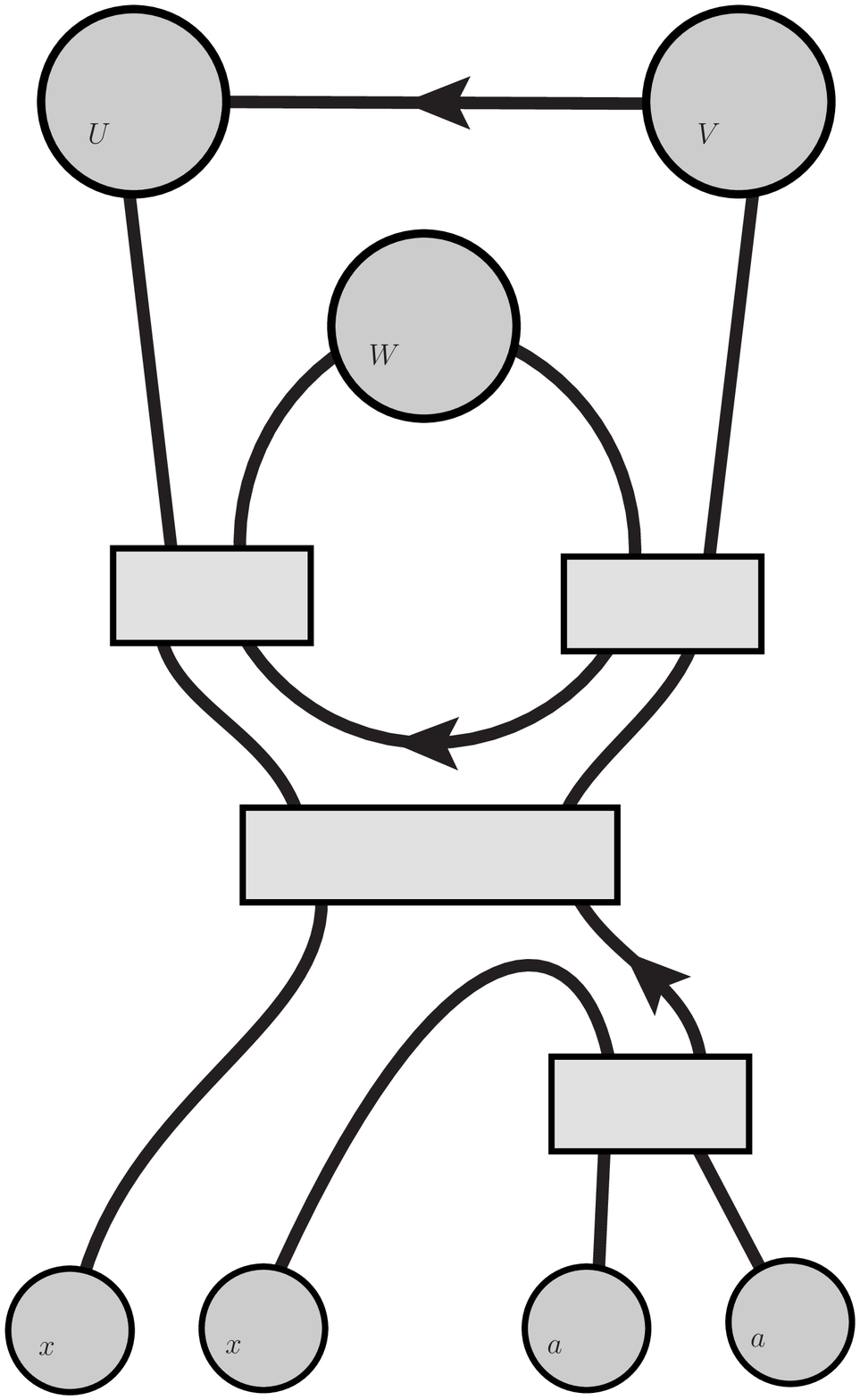}}
\ =-(H_{\rm v},a_{\ux}^2)_1
\ \ 
{\rm with}
\ \ 
H_{\rm v}=\ 
\parbox{3cm}{
\psfrag{U}{$\scriptstyle{U}$}
\psfrag{V}{$\scriptstyle{V}$}
\psfrag{W}{$\scriptstyle{W}$}
\psfrag{x}{$\scriptstyle{x}$}
\includegraphics[width=3cm]{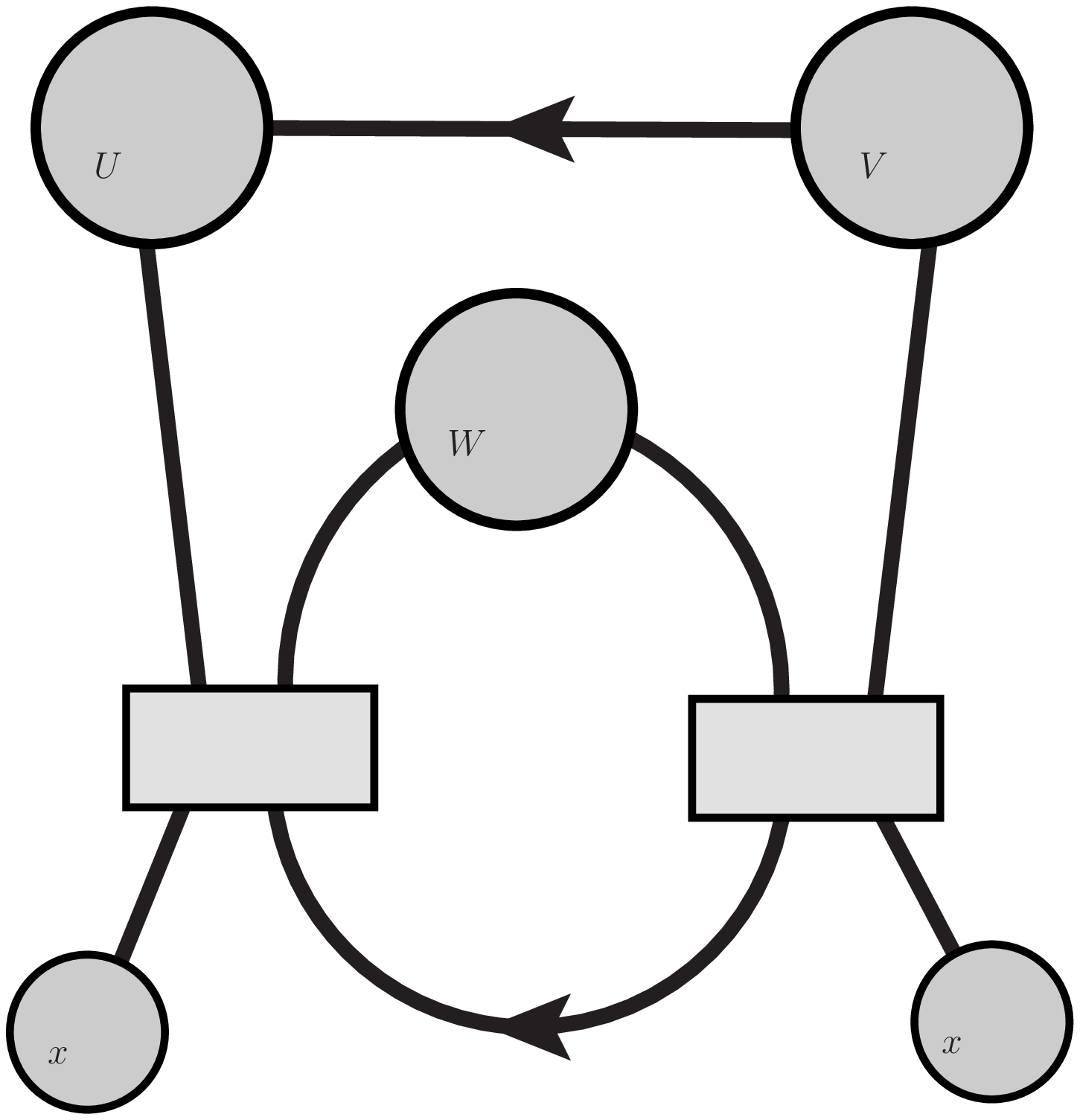}}\ \ .
\]
If we go through the same manipulations for $G_{\rm w}$ and add its
contribution to that of $G_{\rm v}$, then we get the required
expression for $M$ together with 
\[
N=-\frac{1}{2}H_{\rm v}-\frac{1}{2}H_{\rm w}. 
\]
Here $H_{\rm w}$ is the expression similar to $H_{\rm v}$, with $V$
and $W$ interchanged. 
Expanding the two symmetrizers and dropping the zero term with the $W$ self-loop, we get
\[
H_{\rm v}=\frac{1}{4}\ 
\parbox{3cm}{
\psfrag{U}{$\scriptstyle{U}$}
\psfrag{V}{$\scriptstyle{V}$}
\psfrag{W}{$\scriptstyle{W}$}
\psfrag{x}{$\scriptstyle{x}$}
\includegraphics[width=3cm]{pic2.eps}}
\ +\frac{1}{4}
\ 
\parbox{3cm}{
\psfrag{U}{$\scriptstyle{U}$}
\psfrag{V}{$\scriptstyle{V}$}
\psfrag{W}{$\scriptstyle{W}$}
\psfrag{x}{$\scriptstyle{x}$}
\includegraphics[width=3cm]{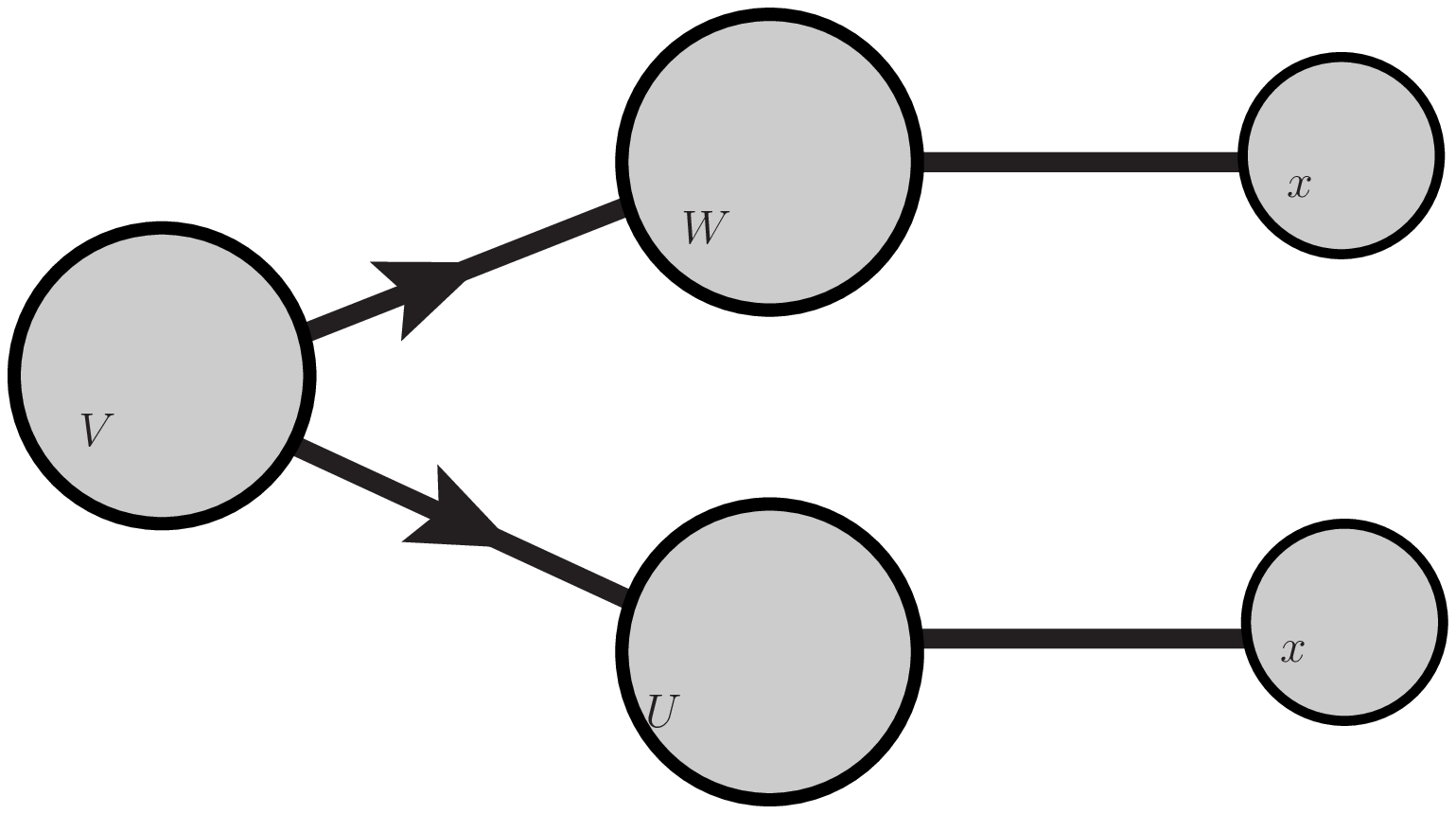}}
\ +\frac{1}{4}\ 
\parbox{3cm}{
\psfrag{U}{$\scriptstyle{U}$}
\psfrag{V}{$\scriptstyle{V}$}
\psfrag{W}{$\scriptstyle{W}$}
\psfrag{x}{$\scriptstyle{x}$}
\includegraphics[width=3cm]{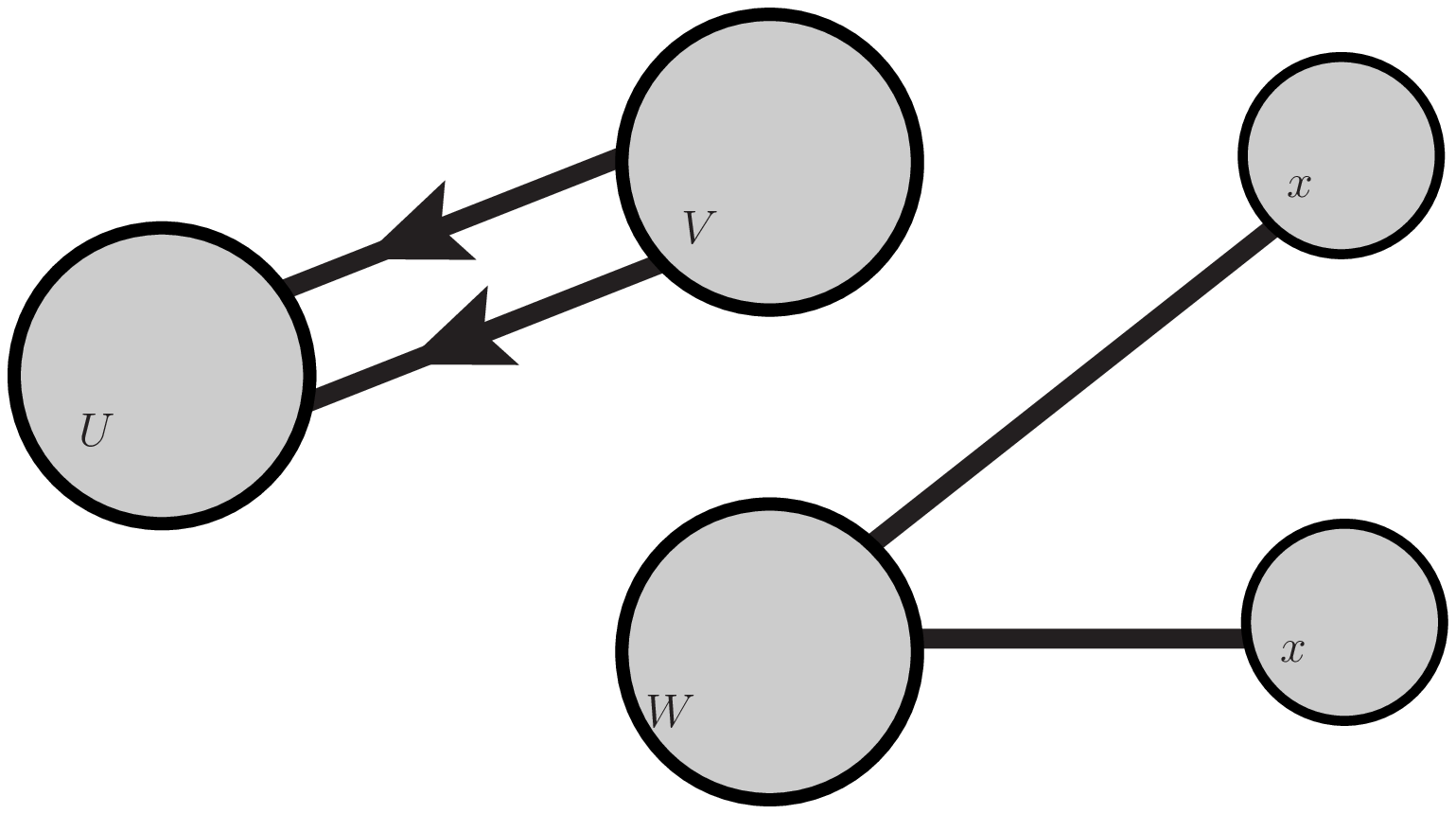}}
\ \ .
\]
By identity (\ref{twochaineq}), the sum of the first two terms is
equal to the last, and thus
\begin{equation}
H_{\rm v}=\frac{1}{2}(U,V)_2 W\ .
\label{Hveq}
\end{equation}
We have seen in (\ref{psiexpeq}) that
\[
(U,VW)_2=\frac{1}{6} (U,V)_2 W+\frac{1}{6} (U,W)_2 V+\frac{2}{3} \{V\rightarrow U\leftarrow W\}\ .
\]
Inserting (\ref{onechaineq}) in the last equation,
we get
\[
(U,VW)_2=\frac{1}{2} (U,V)_2 W+\frac{1}{2} (U,W)_2 V-\frac{1}{3} (V,W)_2 U\ .
\]
By (\ref{Hveq}) and the analogous expression for $H_{\rm w}$, we obtain
\[
(U,VW)_2=-2N-\frac{1}{3} (V,W)_2 U, 
\]
which gives the required expression for $N$. This completes the proof
of Proposition~\ref{proposition.ax}. \qed

\medskip 

\centerline{---} 

\vspace{1cm} 

\parbox{7cm} 
{Abdelmalek Abdesselam \\
Department of Mathematics, \\
University of Virginia, \\
P. O. Box 400137, \\
Charlottesville, VA 22904-4137, \\
USA.\\
{\tt malek@virginia.edu}} 
\hfill 
\parbox{7cm} 
{Jaydeep Chipalkatti \\ 
Department of Mathematics, \\ 
Machray Hall, \\ 
University of Manitoba, \\ 
Winnipeg, MB R3T 2N2, \\ Canada. \\ 
{\tt jaydeep.chipalkatti@umanitoba.ca}}


\bigskip 

\begin{thebibliography}{99} 

\bibitem{Abdesselam}
A.~Abdesselam.
\newblock {On the volume conjecture for classical spin networks.}
\newblock {\sl J. Knot Theory Ramifications}, vol.~21, no.~3, 1250022, 62 pp., 2012. 

\bibitem{Baker} 
H.~F.~Baker. 
\newblock \emph{Principles of Geometry, vol.~II.} 
\newblock Cambridge University Press, 1923. 

\bibitem{Chipalkatti} 
J.~Chipalkatti. 
\newblock {On the coincidences of Pascal lines.} 
\newblock {\sl Forum Geometricorum}, vol.~16, pp.~1--21, 2016. 

\bibitem{Clebsch} 
A.~Clebsch. 
\newblock \emph{Theorie der Bin{\"a}ren Algebraischen Formen. } 
\newblock B.~G.~Teubner, Leipzig, 1872.

\bibitem{ConwayRyba} 
J.~Conway and A.~Ryba. 
\newblock {The Pascal mysticum demystified.} 
\newblock  {\sl Math.~Intelligencer}, vol.~34, no.~3, pp.~4--8, 2012. 

\bibitem{GraceYoung}
J.~H. Grace and A.~Young. 
\newblock \emph{The Algebra of Invariants.}
\newblock Reprinted by Chelsea Publishing Co., New York, 1962. 

\bibitem{Harris} 
J.~Harris. 
\newblock Galois groups of enumerative problems. 
\newblock {\sl Duke Math. J.}, vol.~46, no. 4, pp.~685--724, 1979. 

\bibitem{HowardMSV}
B.~Howard, J.~Millson, A.~Snowden and R.~Vakil.
\newblock The equations for the moduli space of $n$ points on the line.
\newblock {\sl Duke Math. J.}, vol.~146, no.~2, pp.~175--226, 2009. 

\bibitem{KadisonKromann} 
L.~Kadison and M.~T.~Kromann. 
\newblock \emph{Projective Geometry and Modern Algebra.} 
\newblock Birkh{\"a}user, Boston, 1996. 

\bibitem{KungR}
J.~P.~S.~Kung and G.-C.~Rota.
\newblock The invariant theory of binary forms.
\newblock {\sl Bull. American Math. Soc. (N.S.)}, vol.~10, no.~1, pp.~27--85, 1984. 

\bibitem{Olver}
P.~Olver. 
\newblock \emph{ Classical Invariant Theory}. 
\newblock London Mathematical Society Student Texts. Cambridge University 
  Press, 1999.

\bibitem{Pedoe} 
D.~Pedoe. 
\newblock {How many Pascal lines has a sixpoint?}
\newblock {\sl The Mathematical Gazette}, vol.~25, no.~264,
pp.~110--111, 1941. 

\bibitem{PedoeC} 
D.~Pedoe. 
\newblock \emph{Geometry, A Comprehensive Course}. 
\newblock Reprinted by Dover Publications, New York, 1988. 

\bibitem{Salmon}
G.~Salmon. 
\newblock \emph{A Treatise on Conic Sections.}
\newblock Reprint of the 6th ed.~by Chelsea Publishing Co., New York,
2005. 

\bibitem{Schreck_etal} 
P.~Schreck, P.~Mathis, V.~Marinkovi{\'c} and P.~Jani{\v c}i{\'c}. 
\newblock {Wernick's list: a final update.} 
\newblock {\sl Forum Geometricorum}, vol.~16, pp.~69--80, 2016. 

\bibitem{Seidenberg} 
A.~Seidenberg. 
\newblock \emph{Lectures in Projective Geometry}. 
\newblock D.~Van Nostrand Company, New York, 1962. 

\bibitem{Wernick} 
W.~Wernick. 
\newblock Triangle constructions with three located points. 
\newblock {\sl Math.~Mag.}, vol.~55, no.~4, pp.~227--230, 1982. 
\end{thebibliography}
\end{document}